\magnification1200
\def\jump{\bigskip\bigskip}

\def\ra{right-angled }
\def\hn{{\bf H}^n}

\centerline{\bf Boundaries of right-angled hyperbolic buildings
\footnote{$^*$}{\rm Both 
authors were partially supported by KBN grant 
2 P03A 017 25. The second author 
was a Marie Curie Intra-European fellow, contract MEIF CT2005
010050. 
}
}
\smallskip
\centerline{Jan Dymara $^1$, Damian Osajda $^{1,2}$ 
} 
\smallskip
{$^1$ Instytut Matematyczny, Uniwersytet Wroc\l awski, 
pl. Grunwaldzki 2/4, 50--384 Wroc\l aw, Poland}

{$^2$ Institut de Math\'ematiques de Jussieu, Universit\'e Paris 6,
Case 247, 4 Place Jussieu, 75252 Paris Cedex 05, France}
\medskip
\bf Abstract. \rm We prove that the boundary of a right-angled
hyperbolic building is a universal Menger space. Corollary:
the 3-dimensional universal Menger space is the boundary
of some Gromov-hyperbolic group.
\smallskip
\bf Mathematics Subject Classification (2000): \rm 20E42, 54F35, 20F67

\bigskip
\centerline{\bf Introduction}
Hyperbolic \ra buildings were first explored by Marc Bourdon. 
The easiest non-trivial example can be glued from
infinitely many pentagons. We glue them along edges---a finite number greater than 2 along each edge---so 
that a small neighbourhood of each vertex is a cone over a full bipartite graph.
We want the obtained polyhedral complex to be connected
and simply connected; this is easily arranged by passing
to the universal cover of a connected component.
A natural metric on our complex is a piecewise hyperbolic metric,
each pentagon given the shape of a \ra hyperbolic pentagon.
This and similar examples were constructed and
thoroughly investigated by Bourdon ([Bd1], [Bd2])
and Bourdon and Pajot ([BP1], [BP2]).

In particular, Bourdon states that the Gromov boundary 
of any of the complexes he considers is the Menger curve. There are two folklore proofs of this statement.
The first follows the arguments of Benakli (cf. [Bd1], [B]). 
The second uses the result of Kapovich and Kleiner ([KK]): 
if the boundary of a one-ended hyperbolic group is $1$-dimensional and has no local cut
points then it is either the Menger curve or the Sierpi\'nski carpet.
This result applies to uniform lattices in the isometry groups of Bourdon's buildings.
Since (as one can check) the boundary of a thick building 
contains a non-planar subset, it follows that in this case
the boundary is the Menger curve. No details of either of the arguments have been published.

The purpose of this paper is to prove the following theorem.

\smallskip
{\bf Main Theorem}\par\sl
Let $X$ be a locally finite \ra thick hyperbolic building of dimension $n\ge2$.
Then the Gromov boundary of $X$ is homeomorphic to
the universal $(n-1)$-dimensional Menger space $\mu^{n-1}$.
\smallskip\rm

A discussion of buildings explaining the meaning of our assumptions 
is contained in section 1. Let us just mention that
hyperbolic means Lobatchevsky hyperbolic rather than Gromov hyperbolic.
The existence of $X$ as in the theorem is equivalent to the existence
of a bounded finite \ra polyhedron in $\hn$. Therefore, 
$X$ exists only for $n=2,3,4$ ([Vin]). 
Once the polyhedron is given, $X$ can be constructed 
as the universal cover of some finite complex (cf. [D2], [GP] or section 4).
The fundamental group of this finite complex is quasi-isometric to $X$;
thus, one gets examples of Gromov hyperbolic groups with 
Gromov boundary $\mu^3$, $\mu^2$ and $\mu^1$.
The latter two spaces have been known to be boundaries of Gromov
hyperbolic groups (cf. [BK]), but $\mu^3$ is new. 

The proof of the main theorem is based on the following characterisation
of $\mu^{n-1}$, due to Bestvina [Be].
A metric space $Y$ is homeomorphic to $\mu^{n-1}$ if and only if
it is compact, $(n-1)$-dimensional, 
$(n-2)$-connected,
locally $(n-2)$-connected, 
and has the $(n-1)$-dimensional 
disjoint disc property ($DD^{n-1}P$).
We check that these conditions are satisfied for the boundary
$\partial X$ of an $n$-dimensional \ra hyperbolic building $X$.
In Lemma 3.1 we prove that $\partial X$ is compact and $(n-1)$-dimensional.
Recall that $Y$ is locally $(n-2)$-connected
if for each $y\in Y$ and every open neighbourhood $U$ of $y$ 
there exists an open set $V$, $x\in V\subseteq U$,
such that every 
map $S^k\to V$ extends to a 
map $B^{k+1}\to U$
(for $k=0,1,\ldots,n-2$).
In Proposition 3.5 we check that $\partial X$ is 
$(n-2)$-connected and
locally $(n-2)$-connected. 
The $(n-1)$-dimensional disjoint disc property says that for any 
two 
maps $f,g\colon D^{n-1}\to Y$ and any $\epsilon>0$ there exist maps 
$f',g'\colon D^{n-1}\to Y$, $f'$ $\epsilon$-close to $f$,
$g'$ $\epsilon$-close to $g$, such that $f'(D^{n-1})\cap
g'(D^{n-1})=\emptyset$.
A standard way to prove $DD^{n-1}P$ 
is to construct, for any $\epsilon>0$,
two maps $\phi,\psi\colon Y\to Y$, both $\epsilon$-close to the identity map
and with disjoint images. Such maps for $\partial X$ are constructed
in Theorem 4.8 and Corollary 4.10.

In section 4, we investigate the structure of \ra buildings, reproving
some results of Globus ([Gl]) and Haglund and Paulin ([HP]). 
The advantages of our approach are as follows: 
(i) our assumptions on the thickness of the buildings are weaker;
(ii) we obtain a deeper understanding of the automorphism group, allowing
us to prove $DD^{n-1}P$. 

Finally, in the appendix we prove that 
an $n$-dimensional hyperbolic (not necessarily right-angled) 
or Euclidean building
is $(n-2)$-connected at infinity.  
This is a special case of some results of [GP] and [DM]. We briefly criticise the arguments
given in those papers. 

A variant of the main theorem has been proved independently by 
A.Dranishnikov and T.Januszkiewicz. We have not seen the
details of their work. 

We are grateful to Mike Davis, Fr\'ed\'eric Paulin,
 Tadeusz Januszkiewicz and Jacek
\hbox{\'Swi\hskip2.8pt\lower6.1pt\hbox{`}\hskip-5.8pt a\-tkowski} 
for helpful conversations. 

\jump

\def\ra{right-angled }
\def\hn{{{\bf H}^n}}
\def\Z{{\bf Z}}
\def\R{{\bf R}}
\def\bH{{\bf H}^n}
\def\lk{\left\{}
\def\rk{\right\}}

\centerline{\bf 1. Generalities on buildings}
\medskip 

Two standard references for buildings are [Br] and [Ron].
Metrics on buildings are discussed in [D2], and hyperbolic buildings in [GP].

A \it Coxeter system \rm  is a pair $(W,S)$, where $W$ is a group, $S$ is a generating subset of
$W$, and $W=\langle S\mid \{(st)^{m_{st}}\}_{s,t\in S}\rangle$. The numbers $m_{st}$ are positive integers or infinity;
$m_{st}=1$ exactly when $s=t$; $m_{st}=\infty$ means that there is no relation between $s$ and $t$.
We usually speak about a Coxeter group $W$, 
in fact meaning some Coxeter system $(W,S)$.
A Coxeter group $W$ is \it right-angled, \rm if $m_{st}\in\{1,2,\infty\}$ for all $s,t\in S$.
A \it special subgroup \rm of $W$ is a subgroup generated by some subset $T$ of $S$:
$W_T=\langle T\rangle$. It is well known that $(W_T,T)$ is a Coxeter system.
A subset $T\subset S$ is called \it spherical, \rm if $W_T$ is finite.
For example, $\emptyset$ is spherical; $\{s\}$ is always spherical;
$\{s,t\}$ is spherical unless $m_{st}=\infty$.
If $W$ is right-angled, then $T$ is spherical if and only if every two elements of $T$
commute. 
For $w\in W$ we define $\ell(w)$ as the length of a shortest word in the generators $S$
representing $w$. We put $In(w)=\{s\in S\mid \ell(ws)<\ell(w)\}$.
It is well known that $In(w)$ is always spherical.  

Several different descriptions of buildings will be useful for us. We start with a combinatorial one.
Let $W$ be a Coxeter group. We equip $W$ with a family $(\sim_s)_{s\in S}$ 
of equivalence relations, defined as follows: $w\sim_sv\iff w\in\{v,vs\}$.
Suppose that $A$ and $B$ are two sets, each equipped with an $S$-indexed
family of equivalence relations. A map from $A$ to $B$ is a morphism
if it preserves each of the equivalence relations; it is 
an isomorphism if it is a bijective morphism and if its inverse
is also a morphism.
A \it $W$-building \rm is  a set (of chambers) equipped with a family
of equivalence relations $(\sim_s)_{s\in S}$ and with a family of subsets (called apartments)
 isomorphic to $W$, such that:
\item{(B1)} any two chambers are contained in some apartment;
\item{(B2)} {if two chambers $x,y$ are both contained 
in apartments $A,A'$, then 
there exists an isomorphism $A\to A'$ fixing $x$ and $y$.}

\noindent
A building is called \it thick \rm if each equivalence 
class of each relation $\sim_s$
has at least three elements.
We will call a building \it locally finite, \rm if
each equivalence class of each relation $\sim_s$
is finite.
For example, $W$ is a locally finite building, but it is not thick. 
Chambers $x,y$ such that $x\sim_sy$ are called \it $s$-adjacent \rm or simply
{\it adjacent}. A \it gallery \rm in a building $X$ is a sequence
of chambers such that each two consecutive elements are adjacent.
A finite gallery is \it minimal, \rm if there is no shorter gallery
with the same extremities. For a subset $T\subseteq S$ and $x\in X$ we
define the \it residue \rm $Res(x,T)$ as the set of all $y\in X$ such that there 
exists a gallery of the form $x=x_0\sim_{s_1}x_1\sim\ldots\sim_{s_k}x_k=y$,
where $s_1,\ldots,s_k\in T$. For $X=W$ the $T$-residue of $x$ is the left $W_T$-coset containing $x$: $Res(x,T)=xW_T$.
In general, it is well known that $Res(x,T)$ is a $W_T$-building.

The notion of folding map is very important for us.
Let $X$ be a $W$-building.
Pick any chamber $B\in X$. By (B1), for any $x\in X$ there exists an apartment $A$ 
such that $B,x\in A$. Let $\iota_A\colon A\to W$ be the unique isomorphism which sends 
$B$ to $1$. Then, by (B2), $\iota_A(x)$ does not depend on $A$. The formula
$\pi_B(x)=\iota_A(x)$ defines the ($B$-based) 
\it folding map \rm $\pi_B\colon X\to W$.
We often abbreviate $\pi_B$ to $\pi$.
Here is a list of some well-known and useful properties of $\pi$.

\item{(F1)} If $x\in X$ and $t\in In(\pi(x))$, then there exists a unique $x^t\in X$
such that $x\sim_tx^t$ and $\pi(x^t)=\pi(x)t$.
\item{(F2)} The image under $\pi$ of a minimal gallery in $X$ starting at $B$
is a minimal gallery in $W$. Conversely, if $x\in X$ then any minimal gallery 
in $W$ from $B$ to $\pi(x)$ is the image under $\pi$ of a unique minimal gallery from 
$B$ to $x$.
 
\noindent
For $x\in X$ we may define its length (meaning the distance from $B$) 
in terms of the folding: $\ell(x)=\ell(\pi(x))$.
\item{(F3)} For any $x\in X$ and $T\subseteq S$ there exists 
a unique shortest chamber 
$y$ in $Res(x,T)$. Moreover, if $z\in Res(x,T)$, then there exists a minimal gallery
from $B$ to $z$ via $y$. The restriction of $\pi$ to $Res(x,T)$  
composed with left multiplication by $\pi(y)^{-1}$ coincides with the $y$-based 
folding map $\pi_y\colon Res(x,T)\to W_T$.
 
\smallskip

Buildings also have geometric realisations. The most general construction is due to Davis.
Let $D$ be a topological space with a family $(D_s)_{s\in S}$ of subspaces ($D$ is a model for
a chamber; $D_s$ is a model for the intersection of two $s$-adjacent chambers). For $p\in D$
we put $S(p)=\{s\in S\mid p\in D_s\}$. Now for any $W$-building $X$ Davis defines
$X_D=X\times D/\sim$, where $(x,p)\sim(y,q)\iff p=q$ and $x\in Res(y,S(p))$. 
The best choice for $D$ is the \it Davis chamber \rm $K$: it is defined as
the geometric realisation of the poset of all spherical subsets of $S$
(including $\emptyset$);
the subspace $K_s$ is the sub-complex spanned by subsets containing $\{s\}$.
We will denote $X_K$ by $|X|$. For $X=W$ one obtains the 
Davis complex $|W|$.
The geometric 
realisation $|Y|$ of a subset $Y$ of $X$  
is the subset of $|X|$ given by $|Y|=\{[(y,p)]_{\sim}\mid y\in Y, p\in K\}$, 
where $[(y,p)]_\sim$ denotes the 
$\sim$-equivalence class of $(y,p)$.
Apartments in $|X|$ are geometric realisations of apartments 
in $X$. The folding map induces a map 
$|\pi|\colon |X|\to|W|$, which is also called the folding map and 
is usually denoted by $\pi$.
Here are some nice properties of $|X|$:
\item{$\bullet$}
if $X$ is locally finite, then $|X|$ is locally compact;
\item{$\bullet$}
$|X|$ is contractible; 
\item{$\bullet$}
$|X|$ carries a piecewise-Euclidean CAT(0) metric (the Moussong metric).

\smallskip

Let $P$ be a convex polytope in the hyperbolic space $\hn$. Suppose that 
each dihedral angle of $P$ is of the form ${\pi\over k}$, where
the positive integer $k$ may vary from angle to angle.
Then the reflections in codimension-one faces of $P$ generate a Coxeter group $W$; Coxeter groups arising in this way 
will be called \it hyperbolic.
 \rm The group $W$ acts 
on $\hn$ with fundamental domain $P$ (this is a theorem of Poincar\' e). 
The barycentric subdivision of $P$ is isomorphic to the Davis chamber $K$ corresponding
to the group $W$. Using this isomorphism one can define a polyhedral structure and a piecewise hyperbolic metric 
on $|X|$ for any $W$-building $X$. Then each apartment in $|X|$ is isometric to $\hn$, with chambers
corresponding to $W$-translates of $P$ in $\hn$. Moreover, the whole building $|X|$ is $CAT(-1)$ (cf. [D2], [GP]).  
A building $X$ (often meaning the geometric realisation $|X|$, equipped with the $CAT(-1)$ metric 
and the polyhedral structure described above)
corresponding to a hyperbolic Coxeter group will be called a hyperbolic building.
If all dihedral angles of $P$ are $\pi\over2$, then 
$P$, $W$ and $X$ are called \it right-angled.\rm

 The \it Gromov boundary \rm $\partial X$ of a hyperbolic building $X$ 
(or, more generally, of a $CAT(-1)$ space, cf. [BH])
can be defined as the set of geodesic rays $\gamma\colon [0,\infty)\to X$ starting at some fixed point $x_0$. 
The topology on $\partial X$ is defined by the basis of open sets $\{U_r(x)\mid x\in X, r>0\}$, where
$U_r(x)=\{\gamma\in\partial X\mid \gamma([0,\infty))\cap B_r(x)\ne\emptyset\}$
($B_r(x)$ is the open ball in $X$ of radius $r$ centred at $x$).
The topological space thus obtained is independent of the choice of $x_0$.
We will always choose $x_0$ in the interior of a chamber.
For $p,q\in X\cup\partial X$ we denote by $\overline{pq}$ the geodesic
segment from $p$ to $q$ (which exists and is 
unique because $X$ is a $CAT(-1)$ space). 
We define the topology on $X\cup\partial X$ 
by the basis of open sets consisting of
open balls in $X$ and sets $V_r(x)=\{y\in X\cup\partial X\mid \overline{x_0y}\cap B_r(x)\ne \emptyset\}$;
restricted to $\partial  X$ this topology yields the topology described above.
If $X$ is locally compact then $X\cup\partial X$ is a compactification of $X$. 
The folding map $\pi\colon X\to\hn$
extends to a map $\pi\colon X\cup\partial X
\to\hn\cup\partial\hn$, where $\partial{\hn}$ is the Gromov boundary of the hyperbolic space.

\bigskip

\centerline{\bf 2. Half-spaces}
\medskip 
 
The purpose of this section is to prove some auxiliary facts about buildings.
In subsection 2.1 we give a different basis of open sets for the topology 
on $X\cup\partial X$; in 2.2 we prove some properties of the elements
of this new basis; in 2.3 we discuss connectedness properties
of some subsets of spherical buildings. 
\bigbreak

{\bf 2.1. Standard neighbourhoods.}
\medskip

In this subsection we assume that $X$ is a hyperbolic building.
We keep the notation ($W$, $P$, $\pi$, $x_0$, $B_r(x)$,  $V_r(x)$) as in the final two paragraphs of section 1.
In particular, $\pi$ is the folding map based at a chamber $B$ and sending $B$ to $P$. 
We choose $x_0$ in the interior of $B$. 
We denote by $p_R$ the geodesic retraction of $X\cup\partial X$ onto $\overline{B_R(x_0)}$:
$p_R(x)$ is the intersection point of $S_R(x_0)$ and $\overline{x_0x}$ if $d(x_0,x)>R$;
otherwise $p_R(x)=x$. We also use $p_R$ to denote the corresponding retraction in 
$\bH\cup\partial\bH$. 

Let $H$ be a hyperplane in $\hn=|W|$ containing a codimension-one face of some $W$-translate $wP$ of $P$.
Such hyperplane is called a \it wall \rm and divides $\hn$ into two open connected pieces: $H^+$ and $H^-$
(our convention is $int(P)\subseteq H^-$). We put $\partial H^+=\{y\in\partial\hn\mid \overline{py}\cap H^+\ne\emptyset\}$,
for some $p\in int(P)$ (the result does not depend on the choice of such $p$; one may choose $p=\pi(x_0)$).
Let ${\cal H}$ be the set of all connected components of sets of the form
$\pi^{-1}(H^+\cup \partial H^+)$, over all walls $H$. Since all $B_r(x)$ and $V_r(x)$ are pathwise connected 
(any $y\in V_r(x)$ can be connected by a part of 
$\overline{x_0y}$ to a point in $B_r(x)$), the space
$X\cup\partial X$ is locally pathwise connected. 
Therefore, elements of ${\cal H}$ are open; they will be used as neighbourhoods
of boundary points, and called 
{\it standard (open) neighbourhoods}.  
By convention, the whole space $X\cup \partial X$ is also
a standard open neighbourhood. 
\rm We claim that ${\cal H}\cup\{B_r(x)\mid x\in X,r>0\}$ is another basis of open sets
for the topology on $X\cup\partial X$. This is implied by the following lemma.
\smallskip
\bf Lemma 2.1\sl\par
Let $x\in \partial X$ and $U\subset X\cup \partial X$ be its
neighbourhood. Then there exists a wall $H$ 
such that one of connected components
of $\pi ^{-1} (H^+\cup \partial H^+)$ contains $x$ and is contained in
$U$. 
\rm\par
Proof.
By the definition of the topology on $X\cup \partial X$, and because $x_0\in int(B)$, one can find a
point $x_1$ lying on the geodesic ray $\overline {x_0x}$ and inside a
chamber $C$, and a positive number $r$, 
such that
$B_r(x_1)\subseteq int(C)$ and $V_r(x_1)\subseteq U$.
Then $V_r(x_1)$
is an open
neighbourhood of
$x$ contained in $U$. 

{\it Claim:} There exists a wall $F\subset \bH$ such that
$\pi(x)\in
F^+ \cup \partial F^+ \subset \pi (V_r(x_1))$.

{\it Proof of the claim:} 
Notice that $\pi(V_r(x_1))$ equals $V_r(\pi(x_1),\bH)$, 
hence is an open neighbourhood of $\pi(x)$ in $\hn\cup\partial\bH$.
Indeed, this follows easily from the following two observations: 
$\pi(B_r(x_1))=B_r(\pi(x_1),\bH)$;
any geodesic segment in $X$ starting at $x_0$ is mapped by $\pi$
onto a geodesic segment in $\bH$. 

Choose $R>0$ such that every hyperplane
$G\subset \bH$ contained in $\bH \setminus B_R(\pi (x_0),\bH)$ and
intersecting
the geodesic ray $\overline {\pi (x_0)\pi (x)}$ is contained in $\pi(V)$.
Now observe that there
exists a constant $M$ such that if a
hyperplane $G$ satisfies 
$G\cap \overline {\pi(x_1)\pi (x)}\neq \emptyset$ and $G\cap
B_R(\pi (x_0),\bH)\neq \emptyset$, then 
$G\cap B_M(\pi(x_1), \bH)\neq \emptyset$.
There are only finitely many walls in $\hn$ satisfying the latter condition.
Hence, one can choose $R'>R$ such that if $p\in \overline
{\pi (x_0)\pi (x)}$ and $d_{\bH}(\pi (x_0),p)=R'$, then any wall
$F$ intersecting the geodesic ray $\overline {p \pi (x)}$ is
contained in $\bH \setminus B_R(\pi (x_0), \bH)$, hence in $\pi (V_r(x_1))$.
To finish the proof of the claim choose one of such walls.\hfill{$\diamond$(\it Claim)}

Let now $F$ be as in the claim. 
Let $\widetilde F^+$ be the connected component of $\pi^{-1} (F^+\cup\partial F^+)$
containing $x$. 
Let $D$ be the length of a geodesic segment
$\overline{\pi(x_0)q}$ which is tangent to $S_r(\pi(x_1),\bH)$ at $q$.
Then we have $\pi (p_D (\widetilde F^+))\subseteq p_D
(F^+\cup\partial F^+)\subset B_r(\pi(x_1),\bH)$, hence $p_D(\widetilde F^+)\subset \pi^{-1}
(B_r(\pi(x_1),\bH))$. Recall that $B_r(x_1)$ is contained in the interior of one chamber.
Therefore, $B_r(x_1)$ is one of the connected components of $\pi^{-1}(B_r(\pi(x_1),\bH))$.
However, $\widetilde F^+$ is connected and its
closure contains $x$, hence $p_D(\widetilde F^+)$ is connected and contains
$p_D(x)$. Since $p_D(x)=p_D(x_1)$, we have $p_D(\widetilde F^+)\subset B_r(x_1)$. 
This implies  $\widetilde F^+
\subset V_r(x_1)
\subset U$.
\hfill{$\diamond$}

\bf Lemma 2.2\sl\par
Let $H$ be a wall, and let $x\in\pi^{-1}(H^+)$. Then the connected component 
of $x$ in $\pi^{-1}(H^+)$ is dense in the connected component of $x$ in $\pi^{-1}(H^+\cup \partial H^+)$.
\rm \par
Proof. Since the space $X\cup\partial X$ is locally pathwise connected, connected components
of open sets in this space are pathwise connected. Let $y$ be in the connected component
of $x$ in $\pi^{-1}(H^+\cup\partial H^+)$, and let $\gamma\colon [0,1]\to X\cup\partial X$ be a path
from $x$ to $y$ contained in that component. 
We may choose
$R$ so large that $p_R\circ\gamma$ is a path starting at $x$ and contained in $\pi^{-1}(H^+)$.
Concatenating this path with $\overline{p_R(y)y}$ we obtain a path from $x$ to $y$ which is contained
in $\pi^{-1}(H^+)$, except perhaps for its endpoint $y$. It follows that $y$ belongs to the closure of 
the component of $x$ in $\pi^{-1}(H^+)$.
\hfill{$\diamond$}

\bigskip

{\bf 2.2. Shortest elements.}
\medskip
In this subsection $W$ is a \ra Coxeter group (we assume that $W$ is hyperbolic only in Propositions 2.12 and 2.13).
Our goal is to prove that any half-space in $W$ has a unique shortest element (Prop.\thinspace2.5);
we will also investigate the corresponding question for buildings (Prop.\thinspace2.11).
A \it half-space \rm in $W$ is 
a set of the form
$H(w,s)=\{h\in W\mid d(h,ws)<d(h,w)\}$, where $w\in W$, $s\in S$ and $d(w_1,w_2)=\ell(w_1^{-1}w_2)$.
The name is motivated by the fact that if $W$ is hyperbolic then the geometric realisation of $H(w,s)$
is a closed half-space in the usual sense in $|W|=\hn$. 

We begin with two principles which are very useful when dealing with distances in Coxeter groups.
\smallskip
\item{($\pm1$)} {$d(as,b)=d(a,b)\pm1$ and $d(a,bs)=d(a,b)\pm1$, for every $a,b\in W$ and every $s\in S$.}
\smallskip
\item{(R)} {Let $t,t'\in S$ be two distinct commuting generators of $W$, let $R$ be a $\{t,t'\}$-residue in $W$,
and let $x\in W$. Then the four distances from $x$ to chambers of $R$ yield three consecutive
integers, the middle one attained twice, on two non-adjacent chambers of $R$.}
\smallskip
\noindent Property (R) follows from properties ($\pm1$) and (F3) (the latter is stated in section 1). 
Now we proceed to some preliminary lemmas. The proofs are quite standard, so we omit the details.

\smallskip

\bf Lemma 2.3\sl\par
Suppose $h\in H(w,s)\setminus \{ws\}$, $ht\not\in H(w,s)$ for some
$t\in S$. Then there exists $t'\in S$ such that 
$t't=tt'$ and $d(ht',ws)<d(h,ws)$. Moreover:
\item{a)} $ht'\in H(w,s)$, $ht't\not\in H(w,s)$;
\item{b)} $H(ht,t)=H(htt',t)$.
\rm\par
Proof.
Choose $t'$ such that $ht'$ is closer to $ws$ than $h$. 
Using ($\pm1$) one can deduce that $ht'$ is then closer to $w$ than $h$.
Therefore, $t,t'\in In(w^{-1}h)$, hence $tt'=t't$.
It remains to prove a) and b). 

We apply property (R) to the residue $R=Res(h,\{t,t'\})$.
First, we take $x=w$ and $x=ws$. The eight distances are easily determined
up to a common additive constant; part a) follows. Second, take an arbitrary
$x\in W$. Then there are four cases to consider, depending on which
element of $R$ is closest to $x$. In each case it is readily
checked that $x\in H(ht,t)$ if and only if $x\in H(htt',t)$; this proves b).  
\hfill{$\diamond$}
\smallskip

\bf Lemma 2.4 \sl\par
Suppose $h\in H(w,s)$, $ht\not\in H(w,s)$ for some
$t\in S$. Then $t=s$ and $H(hs,s)=H(w,s)$.
\rm\par
Proof. Take a counterexample (to the claim $t=s$) which is closest  to $w$.
Lemma 2.3 produces a counterexample which is even closer 
to $w$, contradiction. Now the same argument proves the second claim.
\hfill{$\diamond$}
\smallskip

\bf Lemma 2.5 \sl\par
Suppose that $h\in H(w,s)$, $hs\not\in H(w,s)$. 
Then there exists a minimal gallery $ws$, $wst_1,\ldots,
wst_1\ldots t_m=h$ such that $t_is=st_i$.
The converse is also true. 
\rm\par
Proof. The first statement follows from Lemma 2.3. 
The converse is easily proved by induction on $m$:
one should apply property (R) to $R=Res(wst_1\ldots t_m,\{t_m,s\})$ 
and $x=w,ws$.  
\hfill{$\diamond$}
\smallskip

\bf Corollary 2.6 \sl\par
a) The set
$\{h\in W\mid h\in H(w,s), hs\not\in H(w,s)\}$ coincides with $wsW_{\{s\}'}$,
where $\{s\}'=\{t\in S\setminus\{s\}\mid ts=st\}$.

b) Any half-space is gallery connected.
\rm\par
Proof. Part a) follows directly from Lemma 2.5.
To prove b), consider a gallery from $x\in H(w,s)$ to $w$. Let 
$y$ be the first element of that gallery which does not belong to $H(w,s)$.
Then, by Lemma 2.4 and part a), $ys\in wsW_{\{s\}'}$, so that it can be 
connected to $ws$ by a 
gallery  in $wsW_{\{s\}'}$. Concatenating the part from $x$ to $ys$ of
the first gallery with the second gallery we obtain a gallery in $H(w,s)$
connecting $x$ to $ws$. Statement b) follows. 
\hfill{$\diamond$}
\smallskip

\bf Proposition 2.7\sl\par
Every half-space (in any \ra Coxeter group) has
a unique shortest element.
\rm\par
Proof.
We assume that $1$ does not belong to our half-space $H(w,s)$---otherwise the statement
is trivial.
Let $x\in H(w,s)$, and let $y$ be the first element in a minimal gallery
from $x$ to $1$ which does not belong to $H(w,s)$. 
Then, by Lemma 2.4 and Corollary 2.6, $ys\in wsW_{\{s\}'}$.
Any residue contains a unique shortest element; let $g$ be the shortest
element in $wsW_{\{s\}'}$. There exists a minimal gallery 
from $ys$ via $g$ to $1$ (cf. property (F3), section 1). Consequently,
$\ell(x)\ge\ell(ys)\ge\ell(g)$; equalities hold only if $x=y=g$.
It follows that $g$ is the unique shortest element in $H(w,s)$.  
\hfill{$\diamond$}
\smallskip

The proof of Proposition 2.7 has the following corollary.
\smallskip
\bf Corollary 2.8\sl\par\nobreak
Any element of $H(w,s)$ can be connected with $1$
by a minimal gallery passing through the shortest element of $H(w,s)$.
\rm\par
Proof. 
In the situation of the proof of Proposition 2.7, concatenate the part from $x$ to $ys$ of the first gallery
with the second gallery. The result is a minimal gallery from $x$ via $g$ to $1$.
\hfill{$\diamond$}
\bigskip

We now turn to buildings.
Let $X$ be a $W$-building, and let $\pi\colon X\to W$ be the $B$-based folding map.
We also fix a half-space $H(w,s)$. We assume that $1\not\in H(w,s)$ (because we are 
eventually interested in standard open neighbourhoods) and that
$ws$ is the shortest element of $H(w,s)$ (we may do so because of Lemma 2.4). 
\smallskip

\bf Lemma 2.9\sl\par
Suppose that $x\in X$, $\pi(x)\in H(w,s)$. 
Then there exists a minimal gallery $(x_0=x,x_1,\ldots, x_\ell=B)$ such that $\pi(x_k)=ws$
for $k=d(\pi(x),ws)$. If $(x_0'=x_0,x'_1,\ldots, x'_\ell)$ is another such gallery, then 
$x'_k=x_k$.
\rm\par
Proof.
Let $\pi(x)=g$. If $\sigma=(g,gs_1,\ldots, gs_1\ldots s_\ell)$ is a minimal
gallery from $g$ to $1$ via $ws=gs_1\ldots s_k$, then 
$\widetilde{\sigma}=(x,x^{s_1},(x^{s_1})^{s_2},\ldots)$ 
is the unique gallery from $x$ to $B$ that folds onto $\sigma$
(cf. properties (F1) and (F2) of the folding map). 
Moreover, $\widetilde{\sigma}$ is minimal.
This gives the first assertion.   

Now let $\tau=(g=\pi(x'_0),\pi(x'_1),\ldots, \pi(x'_\ell)=1)$. By Tits' solution of the word problem in
Coxeter groups, the gallery $\sigma$ 
can be transformed into $\tau$ by a sequence of moves of the form
$$
\matrix{\eta=(\ldots,h_i,h_{i+1}=ht,h_{i+2}=htu,\ldots)\cr
\downarrow\cr
\xi=(\ldots,h'_i=h_i,h'_{i+1}=hu,h'_{i+2}=hut=htu=h_{i+2},\ldots)}\leqno(*)$$
where $tu=ut$. Moreover, this can be done with the $k$'th chamber of the gallery
equal to $ws$ throughout the process (just operate separately on the gallery segments 
from $g$ to $ws$ and from $ws$ to $1$). 
Let $(\sigma=\sigma_1,\ldots,\sigma_m=\tau)$ be a sequence of galleries corresponding to such transformation.
Notice that in the situation of the move $(*)$ the galleries $\widetilde{\eta},\widetilde{\xi}$ 
coincide except for the $(i+1)$-st chamber. This is because both $(y^u)^t$ and $(y^t)^u$ are the shortest
element in $Res(y,\{t,u\})$ so that they coincide (here $y$ denotes the common $i^{\rm th}$ element 
of  $\widetilde{\eta}$ and $\widetilde{\xi}$).
It follows that the $k$'th chamber of each  $\widetilde{\sigma_i}$ is the same; therefore
$x_k=x'_k$.
\hfill{$\diamond$}
\medskip

We now define combinatorial and geometric counterparts of standard open neighbourhoods,
in the setting of general \ra buildings.
Let $x\in X$, $\pi(x)\in H(w,s)$. We define $Y\subseteq X$  
as follows: $y\in Y$ if there exists a gallery $(x=x_0\sim_{s_1}x_1\sim_{s_2}\ldots\sim_{s_m}x_m=y)$ such that
$Res(\pi(x_i),s_i)\subseteq H(w,s)$. We also define $H(w,s)_r=\{[g,p]\in |W|\mid Res(g,S(p))\subseteq H(w,s)\}$
and $Y_r=\pi^{-1}(H(w,s)_r)\cap |Y|=\{[y,p]\in |X|\mid y\in Y,\,Res(\pi(y),S(p))\subseteq H(w,s)\}$.
Notice that $Y=\{y\in X\mid int(|y|)\subseteq Y_r\}$. 
\smallskip

\bf Lemma 2.10\sl\par	
$Y_r$ is pathwise connected and both closed and open in $\pi^{-1}(H(w,s)_r)$. 
\rm \par     
Proof. 

1) $Y_r$ is pathwise connected: 
Recall that the Davis chamber $K$ is the geometric realisation of the poset of all
spherical subsets of $S$. The vertex of $K$ corresponding to $\emptyset$ will be denoted $bar(K)$,
and the vertex corresponding to $\{s\}$ by $bar(K_s)$. The corresponding points in a chamber $|z|\subseteq |X|$
will be denoted $bar(|z|)$, $bar(|z|_s)$. Any point in $|z|$ can be connected to $bar(|z|)$ by a 
line segment contained in $|z|$. 

Let $[y,p]\in Y_r$; the segment
from $[y,p]$ to $bar(|y|)$ is contained in $Y_r$. Now let 
$(x=x_0\sim_{s_1}x_1\sim_{s_2}\ldots\sim_{s_m}x_m=y)$ be a gallery as in the definition of $Y$.
The piecewise linear path $bar(|x_0|)-bar(|x_1|_{s_1})-bar(|x_1|)-\ldots-bar(|y|)$
is contained in $Y_r$: the only problematic points are $p_i=bar(|x_i|_{s_i})$;
however, $S(p_i)=\{s_i\}$ and $Res(\pi(x_i),s_i)\subseteq H(w,s)$.
Thus, any point in $Y_r$ can be connected by a path with $bar(|x|)$ 
so that $Y_r$ is pathwise connected. 

2) $Y_r$ is open: Let $[y,p]\in Y_r$. Let $V$ be the subset of the Davis chamber $K$ consisting of points 
$q$ such that $S(q)\subseteq S(p)$. Then  $\bigcup_{z\in Res(y,S(p))} z\times V$ is an open subset of
$X\times K$, closed under the equivalence relation defining $|X|$. Therefore $\bigcup_{z\in Res(y,S(p))} \{[z,v]\mid v\in V\}$
is an open neighbourhood of $[y,p]$. This neighbourhood is contained in  $Y_r$.
 
3) $Y_r$ is closed in $\pi^{-1}(H(w,s)_r)$: Let $[z,q]$ be in the closure of $Y_r$ in $|X|$. Let $N$ be the open neighbourhood of 
$[z,q]$ constructed in 2); then $N\cap Y_r\ne\emptyset$. Let $[y,p]\in Y_r\cap N$.
Since $[y,p]$ is in the closure of the interior of $|y|$, 
some interior points of $|y|$ belong to $N$. This implies that 
$y\in Res(z, S(q))$, and $[z,q]=[y,q]$. 
We are done, unless $Res(\pi(y),S(q))$ is not contained in $H(w,s)$;
in that case, though, $[\pi(z),q]=[\pi(y),q]\not\in H(w,s)_r$ and
$[z,q]\not\in\pi^{-1}(H(w,s)_r)$.
\hfill{$\diamond$}
\smallskip

\bf Proposition 2.11 \sl\par
$Y$ has a unique shortest chamber.
\rm\par
Proof.
Let $\sigma$ be a minimal gallery from $y\in Y$ to $B$ via
$\pi^{-1}(ws)$; we denote by $a(y)$ the element of $\sigma$ that folds onto 
$ws$ (this element is well defined due to Lemma 2.9).  

Suppose now that $y,y'\in Y$, $y\sim_ty'$ and $Res(\pi(y),t)\subseteq H(w,s)$.
We will prove that $a(y)=a(y')$. If $y'=y^t$, then
there exists a gallery $\sigma$ from $y$ to $B$ via $a(y)$ passing through $y'$
(see the construction of a gallery in the proof of Corollary 2.8: one can start constructing
$\sigma$ by shortening $y$ in an arbitrary manner, provided one does not leave $H(w,s)$),
hence $a(y)=a(y')$ in this case. The case $y=y'^t$ is analogous.
Now suppose that the shortest element $u$ of $Res(y,t)$ 
is different from both $y$ and $y'$. 
Since $y\in Y$ and $Res(\pi(u),t)=Res(\pi(y),t)\subseteq H(w,s)$ we have $u\in Y$.
Then $u=y^t=y'^t$, and
$a(y)=a(u)=a(y')$ by the previous case.
 
It follows that if $y\in Y$ and $(x=x_0,\ldots,x_\ell=y)$ is a gallery as in the definition of $Y$,
then $a(y)=a(x_{\ell-1})=\ldots=a(x)$.
Consequently, any $y\in Y\setminus \{a(x)\}$ is strictly longer than $a(y)=a(x)$.
Thus $a(x)$ is the unique shortest element of $Y$.
\hfill{$\diamond$}
\bigskip
We conclude with two propositions summarising the above discussion
in the hyperbolic case.

\smallskip
\bf Proposition 2.12\sl\par
Suppose that $W$ is a \ra hyperbolic Coxeter group, associated to a polyhedron $P\subseteq \hn$, and let $H$ be a wall.
Then among all $w\in W$ such that $int(wP)\subseteq H^+$ there is a unique shortest one; let us call it $w_0$.
Suppose that $H$ contains the face $w_0P_s$ of $P$; then, for any $x_0\in int(P)$, the geodesic through $x_0$ perpendicular
to $H$ intersects $H$ in the interior of the face $w_0P_s$.
\rm\par
Proof.
Suppose that $H$ contains the face $wP_s$ of the chamber
$wP$. We may assume that $\ell(w)<\ell(ws)$ (swapping $w$ and $ws$ if necessary).
Then, under the usual identification of $|W|$ and $\hn$, the geometric realisation
of $H(w,s)$ corresponds to the closed half-space  $\overline{H^+}$.
This follows easily
from the fact that the distance between two elements of $W$ is equal
to the number of walls separating the corresponding chambers. Now the first assertion of the proposition
follows from Proposition 2.7.

For the second assertion observe that the geodesic $\gamma$ passing through $x_0$ and 
perpendicular to $H$ intersects $H$ at an interior point $\gamma(t)$ of the face $wP_s=wsP_s$. 
Indeed, otherwise $\gamma(t)\in H'$ for some wall $H'\perp H$; then, however,
the image of $\gamma$ is contained in $H'$, and cannot contain $x_0$.
Now if $ws\ne w_0$, then there exists a wall $H'\perp H$ separating $wsP$ from $w_0P$ (and, hence, from $x_0$).
Since $\gamma\perp H$, $\gamma$ does not intersect $H'$. 
On the other hand, $\gamma$ connects points $x_0$ and $\gamma(t)$ lying on different sides of $H'$, a contradiction.
\hfill{$\diamond$}
\medskip 

\bf Proposition 2.13\sl\par
Suppose that $W$ is a \ra hyperbolic Coxeter group, associated to a polyhedron $P\subseteq \hn$, and let $H$ be a wall.
Let $w_0$ be the element of $W$ defined in Proposition 2.12. 
Suppose further that $X$ is a $W$-building, and that $U$ is a connected component of $\pi^{-1}(H^+\cup\partial H^+)$.
Then $\pi^{-1}(w_0P)\cap U$ consists of one chamber.
\rm\par
Proof.
Choose $w\in W$ and $s\in S$ such that $H$ is  the wall separating $wP$ from $wsP$.
We use the notation introduced before Lemma 2.10; we choose $x$ such that $int(|x|)\subseteq U$.
Note that $H(w,s)_r=H^+$. By Lemma 2.10, $Y_r$ is the connected component 
of $\pi^{-1}(H^+)$ that contains $int(|x|)$.
Then Lemma 2.2 implies that $Y_r$ is also equal to the intersection of $U$ and $|X|$.
Recall that $Y=\{y\in X\mid int(|y|)\subseteq Y_r\}$. Therefore, 
the proposition follows from Proposition 2.11. 
\hfill{$\diamond$}
\smallskip

\bigskip

{\bf 2.3. Halves and quarters of spherical buildings.} 
\medskip

In this subsection $Y$ is a finite \ra $W$-building.
Such buildings are spherical, in the following sense. 
Let $\Delta$ be a simplex of dimension $|S|-1$, and let $\Delta_s$
be distinct codimension-one faces of $\Delta$, for $s\in S$ 
(the Davis chamber of $W$ would be isomorphic to a cone over the first barycentric subdivision of $\Delta$).
Then the apartments in $Y_\Delta$ are triangulated spheres.
One equips $W_\Delta=S^{n-1}$ with the standard $CAT(1)$ metric, 
in such a way that each simplex
of the triangulation 
is isometric to a right-angled spherical simplex.
Then one pulls this metric
back by a folding map to a piecewise spherical metric on $Y_\Delta$.
Thus one obtains the standard $CAT(1)$ metric on $Y_\Delta$.
In this subsection we abbreviate $Y_\Delta$ to $Y$.

Buildings as above appear as small spheres around vertices
(or, more generally, as small normal
spheres of cells) in \ra hyperbolic buildings. 
When dealing with complements of balls in a  hyperbolic building or 
in a standard open neighbourhood it is natural to consider 
certain subsets of  spherical buildings. In this subsection we define 
such subsets and prove their higher connectedness in the \ra case.

Let $B\in Y$ be a chamber, 
and let $\pi\colon Y\to S^{n-1}$ be the $B$-based folding map.
We choose $\pi$ so as to have $\pi(B)=\{(x_i)\in S^{n-1}\mid x_1,\ldots,x_n\le0\}$.
The map $\pi$ is then a simplicial map for the following triangulation of $S^{n-1}$:
any simplex $\sigma\subseteq S^{n-1}$ is defined by a conjunction of $n$ conditions of the form
$x_i\le0$, $x_i=0$, $x_i\ge0$, one for each $i$.
Let $C$ be the $(n-1)$-simplex in $S^{n-1}$ which is antipodal to $\pi(B)$, i.e.,
$C=\{(x_i)\in S^{n-1}\mid x_1,\ldots,x_n\ge0\}$. We choose any $v\in int(\pi(B))$; 
in other words, $v$ is a unit vector in $\R^{n+1}$ whose all coordinates are negative. 
We denote by $E^+$ the hemisphere $\{x\in S^{n-1}\mid \langle x,v\rangle\le0\}$, and we put
$Y^+=\pi^{-1}(E^+)$.

\bf Lemma 2.14\sl\par
$\pi^{-1}(C)$ is a deformation retract of $Y^+$.
\rm\smallskip
Proof.
We first construct a deformation retraction $r_t$ from $E^+$ to $C$.
The idea is as follows. If $e\in C$ then put $r_t(e)=e$.
If $e\not\in C\cup\pi(B)$ then there exists a minimal simplex in our
triangulation of $S^{n-1}$ containing $e$; this simplex is a join of some face 
of $\pi(B)$ and some face of $C$. There exists a unique great circle containing $e$ and intersecting 
those two faces; our retraction moves $e$ along that circle towards $C$.
In other words, let $e=e_-+e_+$, where $(e_-)_i=\min(e_i,0)$
and $(e_+)_i=\max(e_i,0)$. We put $r_t(e)=te_-+\sqrt{t^2+{1-t^2\over|e_+|^2}}\,e_+$;
notice that this expression is continuous in  $(t,e)\in[0,1]\times (S^{n-1}\setminus \pi(B))$. 
We have 
$${d\over dt}\langle v,r_t(e)\rangle=
\langle e_-,v\rangle+{t\over\sqrt{t^2+{1-t^2\over|e_+|^2}}}\left(1-{1\over|e_+|^2}\right)\langle e_+,v\rangle.$$ 
This expression is non-negative for $e\in S^{n-1}\setminus \pi(B)$, therefore $e\in E^+$ implies
$r_t(e)\in E^+$.
For an $x\in Y^+$ 
we define $R_t(x)$ as follows: choose any simplex $\Sigma$ of $Y$, containing $x$; 
then $R_t(x)\in\Sigma$, $\pi(R_t(x))=r_t(\pi(x))$.  
Note that if $e\in E^+$ and $e\in\sigma$ for some face $\sigma$ of $S^{n-1}$, then
$r_t(e)\in\sigma$; therefore the definition of $R_t(x)$ does not depend on the choice of
$\Sigma$. 
\hfill{$\diamond$}
\smallskip
It is well known (a self-contained proof is to be found in section 4.1)
that a finite \ra building is a join. More specifically, let
$Y_i=\{x\in Y\mid (\forall j\ne i)(\pi(x)_j=0)\}$. Then 
$Y$ is isomorphic as a simplicial complex and homeomorphic as a topological space
to the join of the sets $Y_i$, $i=1,2,\ldots,n$.
Observe that $\pi^{-1}(C)$ is isomorphic
to the join of the sets $Y_i\setminus B$. It follows that
$\pi^{-1}(C)$ is $(n-2)$-connected (being a join, it is a finite building, 
and thus has the homotopy type of a bouquet of $(n-1)$-spheres).
Now Lemma 2.14 implies the following.
\smallskip

\bf Lemma 2.15\sl\par
$Y^+$ is $(n-2)$-connected.
\rm\smallskip

Now put $Y_i^+=\{x\in Y^+\mid \pi(x)_i\ge 0\}$. 
It is clear that $Y_i^+$ is $R_t$-invariant so that it retracts
to $\pi^{-1}(C)$. In particular, we get as before:
\smallskip  

\bf Lemma 2.16\sl\par
$Y_i^+$ is $(n-2)$-connected.
\rm\par

\jump

\def\q{{\bf q}}

\def\ra{ right-angled }
\def\st{ star-like }
\def\Z{{\bf Z}}
\def\D{{\bf D}}
\def\bH{{\bf H}}
\def\hn{{\bf H}^n}
\def\lk{\left\{}
\def\rk{\right\}}
\def\ov{\overline}

\centerline{\bf 3. Local Connectedness}

\medskip
In this section we prove higher connectedness (local and global) 
of the boundary
of a hyperbolic building.
Let us fix our notation conventions. 
\item{$\bullet$} $W$ denotes a \ra hyperbolic Coxeter group acting on $\hn$
with fundamental domain $P$; 
\item{$\bullet$} $X$ is a locally finite $W$-building
(meaning the Davis realisation with the $CAT(-1)$ metric); 
\item{$\bullet$} $B$ is some fixed chamber of $X$ 
(to be called the base chamber);
\item{$\bullet$} $\pi=\pi_B\colon X\to\hn$ is the $B$-based folding map,
such that $\pi(B)=P$;
\item{$\bullet$} $x_0$ is some fixed generic point in the interior of $B$ (the genericity conditions 
will be specified later);
\item{$\bullet$} $S_R(x,Y)$ and $B_R(x,Y)$ are 
the sphere and the open ball of radius $R$ and centre $x$ in a metric space $Y$.
If $Y=X$, then we use abbreviations $S_R(x)$ and $B_R(x)$. If, additionally,
$x=x_0$, then we write simply $S_R$ and $B_R$.
\item{$\bullet$} $p_R\colon X\cup\partial X\to X$ is the geodesic retraction
onto $\overline{B_R}$, i.e., $p_R(x)$ is the intersection point of
$\overline{xx_0}$ and $S_R$ if 
$d(x,x_0)\ge R$, and $p_R(x)=x$ otherwise. 
\smallskip

\bf Lemma 3.1\sl\par
$\partial X$ is an $(n-1)$-dimensional compactum.
\rm\par
Proof.
Consider an inverse system $\lk (S_k,p_k) \rk _{k=1}^{\infty}$
of spheres centred at $x_0$  with bonding maps
$p_k\colon S_{k+1} \rightarrow S_k$ being the geodesic projections 
onto $S_k$. Then $\partial X=\mathop{\rm inv\,lim} 
\lk (S_k,p_k) \rk$. As every $S_k$ is an $(n-1)$-dimensional
compactum, $\partial X$ is an at most $(n-1)$-dimensional compactum.
But since it contains $S^{n-1}$ (a boundary of an
apartment isometric to $\hn$) it has dimension $n-1$.
\hfill{$\diamond$}

\smallskip
\bf Lemma 3.2\sl\par
Let $U$ be a standard neighbourhood of a point of $\partial X$
and let $R>d(x_0,U)$. Then $\overline U\cap S_R$ is a deformation retract of
$\overline U\setminus B_R$.  
\rm\par
Proof.
Roughly speaking, the retraction is executed by the gradient flow of the restriction 
of the function $d(x_0,\cdot)$ to $\overline{U}\setminus B_R$. The case $U=X$ is easy:
the retraction is $(p_t)_{t\in [R,+\infty]}$, where $p_{+\infty}=id_X$.
 
Let $U$ be a connected component of $\pi^{-1}(H^+\cup\partial H^+)$, for some wall $H\subseteq\hn$. 
We identify $\hn$ with the Poincar\'e disc $\D^n$, in such
a way that $\pi(x_0)$ corresponds to $0$. Then let $Z(x)=-x$ be the vector field on 
$\D^n$ pointing towards $0$. 
We define a vector field $V$ on $\overline{H^+}\setminus B_R(\pi(x_0),\hn)$
as follows. If $x\in H^+\cup\partial H^+\cup\partial H$ then we put $V(x)=Z(x)$.
If $x\in H$ then $V(x)$ is proportional to the orthogonal projection
of $Z(x)$ onto $T_xH$; the proportionality constant is chosen so that  
radial component of $V(x)$ is equal to $Z(x)$. 
The vector field $V$ is not continuous;
nevertheless, it defines a continuous flow $\varphi_V^t$. The trajectory $\varphi_V^t(x)$ 
follows the geodesic $\overline{x\pi(x_0)}$ until it hits $H$; then it moves inside $H$ 
along a geodesic towards the projection of $\pi(x_0)$ onto $H$. The trajectory stops when it
hits $S_R(\pi(x_0),\hn)$ (this may happen before it reaches $H$).
Observe that if a trajectory intersects some wall $H'\ne H$, then
it moves from $(H')^+$ to $(H')^-$. Therefore, the flow $\varphi_V^t$ lifts to a flow
$\psi_V^t$ on $\overline{U}\setminus B_R$. This lift defines a retraction
of $\overline{U}\setminus B_R$ onto $\overline{U}\cap S_R$. 
\hfill{$\diamond$}

\smallskip
\bf Lemma 3.3\sl\par
Let $U$ be a standard neighbourhood of a point of $\partial X$.
Then $\overline U\cap S_{R}$ is $(n-2)$-connected for every $R>0$.
\rm\par
Proof.
Let $U$ be a component of $\pi^{-1}(H^+\cup \partial H^+)$ for some wall $H$.
(The case $U=X\cup\partial X$ is very similar.)
It follows from Propositions 2.12 and 2.13 that for $t$ slightly greater than
$d(\pi(x_0),H)$ the intersection $\overline{U}\cap S_t$ is contained in a single chamber.
This intersection is then a disc, hence is contractible.

Next we would like to understand how the topology of $\overline{U}\cap S_t$ changes 
as $t$ grows. The picture is somewhat reminiscent of Morse theory: the topology 
changes only at some critical radii. 
Suppose that $S_t(\pi(x_0),\hn)$ intersects a (closed) face $\sigma\subseteq H\cup H^+$
of our polyhedral structure at some point $p\in\sigma$. We say that the intersection is critical
if $\sigma$ is perpendicular to $\overline{\pi(x_0)p}$ at $p$; $t$ is then called a critical radius.
We make a generic choice of $x_0$ to ensure that critical intersections
occur only at interior points of the corresponding faces ($p\in int(\sigma)$),
and that to each critical $t$ there corresponds a unique critical intersection.
Notice that $\sigma$ can be a vertex of our polyhedral structure. 
Let $d(\pi(x_0), H)=t_0<t_1<t_2<\ldots$ be the sequence of all critical radii.
It is clear that for $t,t'\in (t_i,t_{i+1})$ the spaces $S_t\cap\overline{U}$
and $S_{t'}\cap\overline{U}$ are homeomorphic (cf. [BMcCM]). We will show 
that for every $i$ and every sufficiently small positive $\epsilon$ 
the space $S_{t_i+\epsilon}\cap\overline{U}$ is $(n-2)$-connected
provided $S_{t_i-\epsilon}\cap\overline{U}$ is $(n-2)$-connected.    

We will first deal with the simplest case: the face $\sigma$ corresponding to 
$t_i$ is a vertex $p\in H^+$. 
Let $Res(p)$ 
be the union of all faces in $\hn$ which contain $p$. We choose $\delta>0$
such that the sphere $D\colon=S_\delta(p,\hn)$ is contained in 
$int(Res(p))$. Let $H_p$ be the hyperplane passing through $p$
and orthogonal to $\overline{\pi(x_0)p}$. This hyperplane divides $D$ into 
two hemispheres, $D^-$ (the one closer to $\pi(x_0)$) and $D^+$. 
There exists an $\epsilon\in(0,\delta)$ such that
$D\cap H_p=D\cap S_{t_i+\epsilon}(\pi(x_0),\hn)$; if necessary, we decrease $\delta$
so as to have $\epsilon<\min\{t_{i+1}-t_i,t_i-t_{i-1}\}$.
The sphere $D$ inherits a triangulation from the polyhedral structure on $\hn$.
We want $S_{t_i+\epsilon}(\pi(x_0),\hn)$ and $S_{t_i-\epsilon}(\pi(x_0),\hn)$
to intersect this triangulation `in the same way'. More precisely, we require 
that there be a homeomorphism of $D$ mapping each simplex into itself and
transforming $D\cap S_{t_i+\epsilon}(\pi(x_0),\hn)$ into $D \cap S_{t_i-\epsilon}(\pi(x_0),\hn)$.
This condition can be achieved by further decreasing $\delta$ (and, consequently, $\epsilon$).

Next we pass to the building. 
By Lemma 3.2, $\overline{U}\cap S_{t_i\pm\epsilon}$ is homotopy 
equivalent to $\overline{U}\setminus B_{t_i\pm\epsilon}$. Let $\pi^{-1}(p)=\{p_1,\ldots,p_k\}$, and let 
$D_j=\pi^{-1}(D)\cap \overline{B_\delta(p_j)}$,   
$D_j^+=\pi^{-1}(D^+)\cap \overline{B_\delta(p_j)}$. 
We have $D_j^+=D_j\setminus B_{t_i+\epsilon}$.   
Put $K_j=\overline{B_\delta(p_j)}$, and let $Y^+$ be the closure of
$(\overline{U}\setminus B_{t_i+\epsilon})\setminus\bigcup_{j=1}^kK_j$.
Furthermore, let $Y_j^+=Y^+\cup K_1\cup\ldots\cup K_j$, for $j=0,1,\ldots,k$.
We will prove, by downward induction on $j$, that $Y_j^+$ is $(n-2)$-connected.
The space $Y_k^+$ is 
homotopy equivalent to
$\overline{U}\setminus B_{t_i-\epsilon}$ (here we need the condition that 
$S_{t_i\pm\epsilon}(\pi(x_0),\hn)$ intersect $D$ `in the same way'), hence it is 
$(n-2)$-connected; 
$(n-2)$-connectedness of $Y^+$ will imply the same property for the homotopy 
equivalent space $\overline{U}\setminus B_{t_i+\epsilon}$.
Observe that the sets $K_j$ are pairwise disjoint, and that $Y_j^+$ is obtained
from $Y_{j-1}^+$ by gluing $K_j$ along $D_j^+$. By Lemma 2.15, $D_j^+$ is 
$(n-2)$-connected, while $K_j$ is clearly contractible. Therefore:
\item{1.} {Connectedness of $Y_j^+$ implies that of $Y_{j-1}^+$.}
\item{2.} {($n>2$) By van Kampen's theorem 
$\pi_1(Y^+_j)=\pi_1(Y^+_{j-1}) 
\ast_{\pi _1(D_j^+)}\pi_1(K_j)$. Since $\pi_1(K_j)=\pi_1(D_j^+)=0$
this implies $\pi_1(Y_{j-1}^+)=\pi_1(Y_j^+)=0$.}
\item{3.} {($n>3$) From the Mayer-Vietoris sequence
$$
\ldots \rightarrow H_{l}(D_j^+)
\rightarrow H_{l}(K_j)
\oplus H_{l}(Y^+_{j-1})
\rightarrow H_{l}(Y^+_j)
\rightarrow H_{l -1}(D_j^+)\rightarrow \dots
$$ 
we get $H_l(Y_{j-1}^+)=H_l(Y_j^+)$ for $l\le n-2$.}

\noindent The conclusion now follows from the Hurewicz theorem.

Now we discuss the general case: $t_i$ is a critical radius, $\sigma$ the corresponding
face, $p$ the intersection point of $S_{t_i}(\pi(x_0),\hn)$ and $\sigma$.
We choose $\delta$ so that $S_\delta(p,\hn)\subseteq int(Res(p))$, and 
we choose $\epsilon$ so that $S_\delta(p,\hn)\cap H_p=
S_\delta(p,\hn) \cap S_{t_i+\epsilon}(\pi(x_0),\hn)$.
Let $\sigma_p^\perp$ be the maximal hyperplane orthogonal to $\sigma$ at $p$.
We put 
$D=\sigma^\perp_p\cap S_\delta(p,\hn)\cap\overline{H^+}$
(intersecting with $\overline{H^+}$ is only necessary if $p\in H$)
and $D^+=D\setminus B_{t_i+\epsilon}(\pi(x_0),\hn)$. 
Again, by decreasing $\delta$ we ensure that 
$\epsilon<\min\{t_{i+1}-t_i,t_i-t_{i-1}\}$ and that
the spheres $S_{t_i\pm\epsilon}(\pi(x_0),\hn)$
intersect $D$ `in the same way'. 
We also set $K=\overline{B_\delta(p,\hn)}\cap\overline{H^+}$  
and $L=\sigma\cap S_\delta(p,\hn)$.

We pass to the building. Let $\pi^{-1}(p)=\{p_1,\ldots,p_k\}$.
We have chosen $\delta$ so small that $\pi^{-1}\bigl(\overline{B_\delta(p,\hn)}\bigr)$
is the disjoint union of $\overline{B_\delta(p_j)}$. We put 
$D_j=\pi^{-1}(D)\cap \overline{B_\delta(p_j)}$,
$D_j^+=\pi^{-1}(D^+)\cap \overline{B_\delta(p_j)}$,
$K_j=\pi^{-1}(K)\cap \overline{B_\delta(p_j)}$,
$L_j=\pi^{-1}(L)\cap \overline{B_\delta(p_j)}$.
Then we define $Y^+$ and $Y^+_j$ exactly as before.
Notice that $K_j$ is homeomorphic to a cone over the join $L_j*D_j$,
and is attached to $Y_{j-1}^+$ along a subset of the base of that cone
homeomorphic to $L_j*D_j^+$. By Lemma 2.15 (if $p\in H^+$) or by Lemma 2.16
(if $p\in H$),  $D_j^+$ is $(n-d-2)$-connected, where $d=\dim(\sigma)$.
Then $L_j*D_j^+$ is $(n-2)$-connected by the suspension theorem
($L_j$ is a $(d-1)$-dimensional sphere). Moreover, $K_j$ is 
contractible. In the remaining part of the argument (1.--3.) we just 
replace $D_j^+$ by (a homeomorphic copy of) $L_j*D_j^+$.
\hfill{$\diamond$}

\smallskip

\bf Lemma 3.4\sl\par
Let $U$ be a standard neighbourhood of $x\in \partial X$
in $X\cup \partial X$. Then for every
standard neighbourhood $V$ of $x$ 
whose closure is contained in $U$, 
every $k\in \lk 0,1,2,...,n-2 \rk$ and every map $f\colon S^{k}\rightarrow
V\cap \partial X$ there exists an extension $g\colon B^{k+1}
\rightarrow
U\cap \partial X$ of $f$.
\rm\par
Proof.
We will proceed by induction on $k$.

1.{\bf k=0}. Let $V$ and $f\colon \lk 0,1 \rk \rightarrow V\cup \partial
X$ be given.
We will construct the desired $g\colon I=[0,1] \rightarrow
U\cap \partial X$ as a limit of a sequence $(g_i)_{i=0}
^{\infty}$ of maps $g_i\colon I\rightarrow S_{N_i}\cap U$,
for an increasing sequence of integers $N_i$.

Let us assume that we defined a natural number $N_i$,
a map $g_i\colon I\rightarrow S_{N_i}\cap U$,
and additionally finite families
$V_i$ and $U_i$ of standard open neighbourhoods 
of points of $\partial X$,
and a triangulation
${\cal T}_i$ of $I$ of mesh at most $2^{-i}$ together with 
a map $h_i\colon |{\cal T}_i^{(0)}|\rightarrow U\cap \partial X$
and a map $s_i\colon {\cal T}_i^{(1)}\rightarrow U_i$. (Note
that by ${\cal T}^{(j)}$ we denote the set of $j$-simplices
of a triangulation ${\cal T}$, and $|{\cal T}^{(j)}|$ denotes
a geometric realization of the $j$-skeleton of ${\cal T}$.)
Assume that they satisfy the following conditions:
\item{i)} $h_i|_{\lk 0,1 \rk}=f$,
\item{ii)} $g_i|_{|{\cal T}_i^{(0)}|}=p_{N_i}\circ h_i$,
\item{iii)} $\forall (\tau \in {\cal T}_i^{(1)})
\exists (B\in V_i) (\overline B\subset s_i(\tau) \; \;
{\rm and} \; \;
g_i(\partial \tau)\subset B\; \; 
{\rm and}\; \;g_i(\tau)\subset \overline B$),
\item{iv)} $\forall (B\in V_i)\; \; \overline{B}\subset U$. 

We will show how to find a natural $N_{i+1}$, a map $g_{i+1}$
etc.
For every $D\in V_i$ one can find 
finite families $U_{i+1}^D$ and $V_{i+1}^D$
of standard open neighbourhoods of points of $\partial X$
such that  
\item{a)} $\overline D\cap \partial X\subset 
\bigcup V_{i+1}^{D}\,$,
\item{b)} $\forall (C\in U_{i+1}^{D}) \; \; \forall (A\in U_i)
\; \; {\rm if} \; \; \overline D\subset A \; \; {\rm then} \; \;
\overline C\subset A$,
\item{c)} $\forall (A\in U_{i+1}^{D})\; \; \overline {A}\subset
U\setminus B_{N_i}$,
\item{d)} $\forall (B\in V_{i+1}^{D})\; \; 
\exists (A\in U_{i+1}^{D})
\; \; \overline B\subset A$.

\noindent Define finite families $V_{i+1}$ and $U_{i+1}$ by
$V_{i+1}=\bigcup_{D\in V_i}V_{i+1}^{D}$ and
$U_{i+1}=\bigcup_{D\in V_i}U_{i+1}^{D}$. 
Find a natural $N_{i+1}>N_i$ such that for every $D\in V_i$
we have $S_{N_{i+1}}\cap \overline D \subset
\bigcup V_{i+1}^{D}$.
Given a $1$-simplex $\tau$ of 
${\cal T}_i$, by iii) we find $D_{\tau}\in V_i$ with 
$\overline {D_{\tau}}\subset s_i(\tau)$, $g_i(\partial \tau)
\subset D_{\tau}$ and $g_i(\tau)\subset \overline {D_{\tau}}$.
Every standard open neighbourhood $D$ has the following property:
for any $R>0$ and any $y\in X\cup\partial X$, if $p_R(y)\in D$ then
$y\in D$. Observe that $p_{N_i}(p_{N_{i+1}}\circ h_i(\partial\tau))=
g_i(\partial\tau)\subset D_\tau$; therefore
$p_{N_{i+1}}\circ h_i|_{\partial\tau}$ maps $\partial\tau$ 
into $S_{N_{i+1}}\cap D_\tau$.
By Lemma 3.3 we can extend this map to 
$g_{i+1}^\tau\colon\tau\to S_{N_{i+1}}\cap D_\tau$.
Define $g_{i+1}$ as the union of $g_{i+1}^\tau$ over all $\tau\in{\cal T}_i^{(1)}$.
By continuity of $g_{i+1}$ one can choose a
subdivision ${\cal T}_{i+1}$ of the triangulation ${\cal T}_i$ 
of $I$ with simplices of diameter
at most $2^{-i-1}$, so fine that for every $1$-simplex
$\sigma$ of ${\cal T}_{i+1}$
contained in a 1-simplex $\tau$ of ${\cal T}_i$ there exists 
$B\in V_{i+1}^{D_{\tau}}$ 
such that $g_{i+1}^{\tau}(\sigma)\subset B$.
Then, by d), for any $\sigma$, $\tau$ and $B$ as in the previous sentence
there exists an $s_{i+1}(\sigma)\in U_{i+1}^{D_\tau}$
satisfying $\overline{B}\subset s_{i+1}(\sigma)$.
Observe that, by b), $\overline{s_{i+1}(\sigma)}\subset s_i(\tau)$. 
Finally, we define $h_{i+1}$ as follows: for $v\in{\cal T}_i^{(0)}$ we put $h_{i+1}(v)=h_i(v)$; 
for $v\in{\cal T}_{i+1}^{(0)}\setminus {\cal T}_i^{(0)}$
we choose any point $h_{i+1}(v)\in\partial X$ such that 
$p_{N_{i+1}}(h_{i+1}(v))=g_{i+1}(v)$.

To start the construction of $g_i$'s one has to set:
$N_0$, $g_0$, $V_0$ and $U_0$, ${\cal T}_0$,
$h_0$ and $s_0$.
Let $N_0$ be a natural number such that $S_{N_0}\cap
V\supset p_{N_0}\circ f (S^0)$. By Lemma 3.3 one can find
a map $g_0\colon I\rightarrow S_{N_0} \cap \overline V$ extending
the map $p_{N_0}\circ f\colon  \lk 0,1 \rk\rightarrow S_{N_0}$. Then set
$V_0=\lk V \rk$, $U_0=\lk U \rk$, ${\cal T}_0$ - triangulation
of $B^1$ consisting of one $1$-simplex,
$h_0=f$ and $s_0(v)=U$ for every $v\in {\cal T}_0^{(0)}$.
Then it is obvious that conditions i) - iv) are satisfied.

We will now show some properties of the sequence 
$(g_i)_{i=1}^{\infty}$ of maps that will imply that its
limit is a continuous map extending $f$.

{\it Claim 0: $g_i(\tau)\subset s_L(\tau)$ for $\tau \in {\cal
T}_L^{(1)}$ and $i\geq L$.} First we show that for $i,j=0,1,2,...$
and for any two simplices $\sigma \in {\cal T}_i^{(1)}$ and
$\rho \in {\cal T}_{i+j}^{(1)}$ such that $\rho \subset
\sigma$ we have
$\overline {s_{i+j}(\rho)}\subset s_i(\sigma)$. We proceed by induction
on $j$. For $j=0$ the inclusion is obvious, and for $j=1$ it
follows from the construction of $s_i$. Assume we have proved
that $\overline {s_{i+j}(\rho)}\subset s_i(\sigma)$. Let
$\kappa \in {\cal T}_{i+j+1}^{(1)}$ be a simplex contained in
a simplex $\rho \in {\cal T}_{i+j}^{(1)}$ that is itself
contained in $\sigma \in {\cal T}_i^{(1)}$. Then, by the induction
assumptions, we have
$
\overline {s_{i+j+1}(\kappa)}\subset s_{i+j}(\rho)\subset s_i(\sigma)
$. This finishes the induction.

Let $A=\{\sigma\in {\cal T}_i^{(1)}\mid \sigma\subset\tau\}$.
Then $\tau = \bigcup_{\sigma
\in A} \sigma$, and
$$
g_i(\tau)=g_i(\bigcup_{\sigma
\in A} \sigma)=\bigcup_{\sigma
\in A} g_i(\sigma)\subset \bigcup_{\sigma
\in A} s_i(\sigma)\subset
\bigcup_{\sigma
\in A} \overline{s_i(\sigma)}\subset
s_L(\tau).
$$
Here the
last inclusion follows from what we proved above and the
first one holds by iii).

{\it Claim 1: For every $y\in I$ the limit 
${\rm lim}_{i\rightarrow \infty}g_i(y)$ exists.} 
Take an arbitrary open (in $X\cup \partial X$) finite cover 
${\cal W}$ of $\overline {U}\cap
\partial X$. For every $j\ge i>0$
and every $A\in U_j$ we have
$A\subset U\setminus B_{N_{i-1}}$. 
Therefore there exists a natural
$L>N_0$ such that for every 
$i\geq L$
every neighbourhood $A\in U_i$
is contained in 
some member of ${\cal W}$.
Take an arbitrary $y\in I$. Let $\tau \in {\cal T}_L$
be a maximal simplex containing $y$. 
Then, by Claim 0, $g_i(\tau)\subset s_L(\tau)\subset
W$ for every $i\geq L$ and some $W\in {\cal W}$. This implies
the existence of the limit.

{\it Claim 2: ${\rm lim}_{i\rightarrow \infty}
g_i(y)\in U\cap \partial X$.} Follows from: $g_i(y)
\in \bigcup U_1$ for every $i$;  
$\overline A\subset U$ for every $A\in U_1$.

{\it Claim 3: The formula
$g(x)={\rm lim}_{i\rightarrow \infty}g_i(x)$ 
defines a continuous map $g\colon I\rightarrow U\cap \partial X$.}
As in the proof of Claim 1, 
for every finite open cover ${\cal W}$ there exists $L>0$ such that
for every $i\geq L$ and any $A\in U_i$ the star
$\bigcup {\rm St}(A)$ of $A$ in $U_i$ is contained in some 
member of ${\cal W}$. 
Take an arbitrary $y\in I$. Let $\tau \in {\cal T}_L$
be a maximal simplex containing $y$. 
As in Claim 1, we have $g_i(\sigma)\subset s_L(\sigma)$ for every 
$i\geq L$ and every $1$-simplex $\sigma$ of ${\cal T}_L$
which has non-empty intersection with $\tau$;
hence, $g_i(\bigcup St(\tau))\subset \bigcup 
\lk s_L(\sigma) | \sigma \in St(\tau)\rk \subset \bigcup St(s_L(\tau))
\subset W$ for some $W\in {\cal W}$.
In other words, for every open cover ${\cal W}$ as above
and any given $y\in I$
there exists a natural $L$, $W\in {\cal W}$, and an open neighbourhood
$E \subset I$ of $y$ such that for every
$i\geq L$ we have $g_i(E)\subset W$. This implies that the limit
of $g_i$'s is continuous.

{\it Claim 4: The map $g\colon I\rightarrow U\cap \partial X$
extends $f$.} Follows from the fact that $g_i|_{\lk 0,1 \rk }=
p_{N_i}\circ h_i|_{\{0,1\}}=p_{N_i}\circ f$
and ${\rm lim}_{i\rightarrow \infty}p_{N_i}\circ f(y)=f(y)$ for
every $y\in \lk 0,1 \rk $.
\smallskip

2.{\bf Induction step.} Assume we have proved the Lemma for 
$k=0,1,...,M-1$. Let $V$ and $f\colon S^M\rightarrow V\cup \partial
X$ be given. 
Again, we will construct the desired $g\colon B^{M+1} \rightarrow
U\cap \partial X$ as a limit of a sequence $(g_i)_{i=0}
^{\infty}$ of maps $g_i\colon B^{M+1}\rightarrow S_{N_i}\cap U$,
where $N_i$ is an increasing sequence of integers.

Let us assume that we defined a natural number $N_i$,
a map $g_i\colon B^{M+1}\rightarrow S_{N_i}\cap U$,
and additionally, for every $p=1,2,...,M+1$, finite families
$V_i^p$ and $U_i^p$ of standard open 
neighbourhoods of points of $\partial X$, 
and a triangulation
${\cal T}_i$ of $B^{M+1}$ of mesh at most $2^{-i}$ together with 
a map $h_i\colon |{\cal T}_i^{(M)}|\rightarrow U\cap \partial X$
and a map $s_i\colon {\cal T}_i^{(M+1)}\rightarrow U_i^1$.
Assume that they satisfy the following conditions:
\item{i)} $h_i|_{S^M}=f$,
\item{ii)} $g_i|_{|{\cal T}_i^{(M)}|}=p_{N_i}\circ h_i$,
\item{iii)} $\forall (\tau \in {\cal T}_i^{(M+1)})
\exists (B\in V_i^1) (\overline B\subset s_i(\tau) \; \;
{\rm and} \; \; g_i(\partial \tau)\subset B\; \; 
{\rm and}\; \;g_i(\tau)\subset \overline B$),
\item{iv)} $\forall (B\in V_i^1) \; \; \overline{B}\subset U$. 

We will show how to find a natural $N_{i+1}$, a map $g_{i+1}$
etc.
For every $D\in V_i^1$ and every $p=1,2,...,M+1$ one can find 
finite families $U_{i+1}^{D,p}$ and $V_{i+1}^{D,p}$
of standard open neighbourhoods of points of $\partial X$
such that  
\item{a)} $\overline D\cap \partial X\subset 
\bigcup V_{i+1}^{D,p}\,$,
\item{b)} $\forall (C\in U_{i+1}^{D,p}) \; \; \forall (A\in U_i^1)
\; \; {\rm if} \; \; \overline D\subset A \; \; {\rm then} \; \;
\overline C\subset A$,
\item{c)} $\forall (A\in U_{i+1}^{D,1})\; \; \overline {A}\subset
U\setminus B_{N_i}$,
\item{d)} $\forall (B\in V_{i+1}^{D,p})\; \; 
\exists (A\in U_{i+1}^{D,p})
\; \; \overline B\subset A$,
\item{e)} $\forall (p\geq 2)\; \; \forall (A\in U_{i+1}^{D,p}) 
\; \;\exists (C\in V_{i+1}^{D,p-1})
\; \; \overline{\bigcup St(A,U_{i+1}^{D,p})}\subset C$.

Define finite families $V_{i+1}^p$ and $U_{r+1}^p$ by
$V_{r+1}^p=\bigcup_{D\in V_r^1}V_{r+1}^{D,p}$ and
$U_{r+1}^p=\bigcup_{D\in V_r^1}U_{r+1}^{D,p}$.
Find a natural $N_{i+1}'>N_i$ such that for every $D\in V_i^1$
we have
$S_{N_{i+1}'}\cap \overline D \subset
\bigcup V_{i+1}^{D,M+1}$. Given an $(M+1)$-simplex $\tau$ of 
${\cal T}_i$, by iii) we find $D_{\tau} \in V_i^1$ with 
$\overline {D_{\tau}}\subset s_i(\tau)$,
$g_i(\partial \tau)
\subset D_{\tau}$ and 
$g_i(\tau)
\subset \overline {D_{\tau}}$. 
Observe that then $h_i(\partial \tau)
\subset D_{\tau} \cap \partial X$ and that, by b), 
$\overline B\subset s_i(\tau)$ for every $B\in V_{i+1}^{D_{\tau},p}$,
$p=1,2,...,M+1$. 
Using Lemma 3.3 
one can find a map 
$g_{i+1}^{'\tau}\colon \tau\rightarrow S_{N_{i+1}'}\cap \overline D_{\tau}$ 
extending
the map $p_{N_{i+1}'}\circ h_i|_{\partial \tau}\colon \partial \tau
\rightarrow
S_{N_{i+1}'}\cap D_{\tau}$. 
By continuity of (every) $g_{i+1}^{'\tau}$, one can choose a
subdivision ${\cal T}_{i+1}$ of the triangulation ${\cal T}_i$ 
of $B^{M+1}$ with simplices of diameter
at most $2^{-i-1}$, so fine that for every $1$-simplex
$\sigma$ of ${\cal T}_{i+1}$
contained in suitable $\tau$ there exists $B\in V_{i+1}^{D_{\tau},M+1}$ 
such that $g_{i+1}^{'\tau}(\partial \sigma)\subset B$ . For every vertex
$v$ of ${\cal T}_{i+1}$ not belonging to $|{\cal T}_i^{(M)}|$ one can 
choose a point $\tilde v\in \partial X$ such that
$p_{N_{i+1}'}(\tilde v)=g_{i+1}'(v)$, where $g_{i+1}'$ is 
the union of maps of the form $g_{i+1}^{'\tau}$ over all maximal
simplices $\tau$ of ${\cal T}_i$. 
For a vertex $v\in |{\cal T}_i^{(M)}|$ we put  $\tilde v=h_i(v)$.
Again, by induction assumptions,
for every two vertices $v,w$ of ${\cal T}_{i+1}$ joined by an edge $\langle v,w\rangle$
contained in $\tau$ and not in $|{\cal T}_i^{(M)}|$, 
and for the corresponding points $\tilde v,
\tilde w\in \partial X$, 
one can find $A\in U_{i+1}^{D_{\tau},M+1}$ and a map
$q\colon \langle v,w\rangle \rightarrow A\cap \partial X$ such that $q(v)=\tilde v
,q(w)=\tilde w$.
Assume we proved that for any $l$-simplex $\sigma$ of ${\cal T}_{i+1}$
contained in $\tau$ and not in $|{\cal T}_{i}^{(M)}|$
there exist $A_0,A_1,...,A_l\in U_{i+1}^{D_{\tau},M+3-l}$ and maps
$q_0,q_1,...,q_l\colon \partial \sigma \rightarrow (A_0\cup A_1\cup ...
\cup A_l)\cap \partial X$ 
sending $(l-1)$-faces of $\sigma$ into distinct $A_i$'s and coherent
on their intersections
(we have just checked this for $l=2$).
Since $\partial \sigma \subset \bigcup St(\kappa)$ for every $(l-1)$-simplex
$\kappa$ of $\partial \sigma$, we have $\bigcup_{i=0}^l A_i
\subset \bigcup St(A_0)$. Thus there exists $B\in V_{i+1}^{D_{\tau},M+2-l}$
such that $(\bigcup_{i=0}^lq_i)(\partial \sigma)\subset B$.
Hence, if $l\leq M$, by induction assumptions there exists
$A\in U_{i+1}^{D_{\tau},M+2-l}$ and
a map $q\colon \sigma \rightarrow A\cap \partial X$
extending $\bigcup_{i=0}^lq_i$.
If $l=M+1$, we conclude that for every $(M+1)$-simplex $\sigma$
of ${\cal T}_{i+1}$ contained in $\tau$ there exists 
$B\in V_{i+1}^{D_{\tau},1}$ and a map
$q\colon \partial \sigma \rightarrow B\cap \partial X$ such that
$q(v)=\tilde v$ for every vertex $v$ of $\sigma$ and $q$ coincides with $h_i$
on $\partial \tau\cap\partial\sigma$. By d), there exists $A\in U_{i+1}^{D_{\tau},1}$
such that $\overline B\subset A$. Define $s_{i+1}\colon {\cal T}_{i+1}^{(M+1)}
\rightarrow U_{i+1}^1$ setting $s_{i+1}(\sigma)=A$.
Observe that since $\overline {D_{\tau}}\subset s_i(\tau)$ by b)
we have $\overline A\subset s_i(\tau)$. In other words
for every $(M+1)$-simplex $\tau$ of ${\cal T}_i$ and an $(M+1)$-simplex
$\sigma \subset \tau$ of ${\cal T}_{i+1}$ we have
$\overline {s_{i+1}(\sigma)}\subset s_i(\tau)$.
Because maps of the form $q$ by definition coincide 
on intersections of their domains, their union is a well-defined continuous map
$h_{i+1}^{\tau}\colon |{\cal T}_{i+1}^{(M)}|\cap \tau \rightarrow
A\cap \partial X$. Note that $h_{i+1}^{\tau}|_{\partial \tau}=
h_i|_{\partial \tau}$, and that for every
$(M+1)$-simplex $\sigma\subset\tau$ of ${\cal T}_{i+1}$ there exists 
$B\in V_{i+1}^{D_{\tau},1}$ satisfying
$h_{i+1}^{\tau}(\partial \sigma)\subset B\cap \partial X$
and $\overline B\subset s_{i+1}(\sigma)$.
Because maps of the form $h_{i+1}^{\tau}$ for different choices
of $\tau$ coincide on intersections of their domains,
we can define $h_{i+1}\colon |{\cal T}_{i+1}^{(M)}|\rightarrow U\cap
\partial X$ as the union of all those maps.
One can find a natural $N_{i+1}>N_{i+1}'$ such that 
for every $\tau$ 
and every
$(M+1)$-simplex $\sigma$ of ${\cal T}_{i+1}$ contained in $\tau$
there exists $B\in V_{i+1}^{D_{\tau},1}$ with
$p_{N_{i+1}}\circ h_{i+1}(\partial \sigma)\subset B\cap 
S_{N_{i+1}}$ and $\overline B\subset s_{i+1}(\sigma)$.
By Lemma 3.3, for every such $\sigma$ and $B$ there exists
a map $g_1^{\sigma}\colon \sigma \rightarrow S_{N_{i+1}}\cap \overline B$
extending the map $p_{N_{i+1}}\circ h_{i+1}|_{\partial \sigma}\colon 
\partial \sigma \rightarrow
S_{N_{i+1}}\cap B$. The union of such maps over all maximal
simplices 
defines 
a map $g_{i+1}^{\tau}\colon \tau \rightarrow S_{N_{i+1}}\cap C_{\tau}$,
which extends the map $p_{N_{i+1}}\circ h_{i+1}\colon \partial \tau 
\rightarrow
S_{N_{i+1}}\cap C_{\tau}$ and hence also the map
$p_{N_1}\circ h_i|_{\partial \tau}\colon \partial \tau \rightarrow
S_{N_{i+1}}\cap C_{\tau}$. We define the map $g_{i+1}\colon B^{M+1} \rightarrow
S_{N_{i+1}}\cap U$ as the union of maps $g_{i+1}^{\tau}$ over all
maximal simplexes $\tau$ of ${\cal T}_i$.
Observe that by the construction $h_{i+1},g_{i+1},{\cal T}_{i+1},V_{i+1}^1,
U_{i+1}^1$ satisfy induction assumptions i)-iv) so that one can 
proceed with following steps of the construction.

To start the construction of $g_i$'s one has to set:
$N_0$, $g_0$, $V_0^1$ and $U_0^1$, ${\cal T}_0$,
$h_0$ and $s_0$.
Let $N_0$ be a natural number such that $S_{N_0}\cap
V\supset p_{N_0}\circ f (S^M)$. By Lemma 3.3 one can find
a map $g_0\colon B^M\rightarrow S_{N_0} \cap \overline V$ extending
the map $p_{N_0}\circ f\colon  S^M\rightarrow S_{N_0}$. Then set
$V_0^1=\lk V \rk$, $U_0^1=\lk U \rk$, ${\cal T}_0$ - triangulation
of $B^{M+1}$ consisting of one $(M+1)$-simplex,
$h_0=f$ and $s_0(\sigma)=U$ for every $\sigma \in {\cal T}_0^{(M)}$.
It is obvious that conditions i) - iv) are satisfied.

We will now show some properties of the sequence 
$(g_i)_{i=1}^{\infty}$ of maps which will imply that its
limit is a continuous map extending $f$.

{\it Claim 0: $g_i(\tau)\subset s_L(\tau)$ for $\tau \in {\cal
T}_L^{(M+1)}$ and $i\geq L$.} First we show that for $i,j=0,1,2,...$
and for any two simplices $\sigma \in {\cal T}_i^{(M+1)}$ and
$\rho \in {\cal T}_{i+j}^{(M+1)}$ such that $\rho \subset
\sigma$ we have
$\overline {s_{i+j}(\rho)}\subset s_i(\sigma)$. We proceed by induction
on $j$. For $j=0$ the inclusion is obvious and for $j=1$ it
follows from the construction of $s_i$. Assume we have proved
that $\overline {s_{i+j}(\rho)}\subset s_i(\sigma)$. Let
$\kappa \in {\cal T}_{i+j+1}^{(M+1)}$ be a simplex contained in
a simplex $\rho \in {\cal T}_{i+j}^{(M+1)}$ that is itself
contained in $\sigma \in {\cal T}_i^{(M+1)}$. Then, by the induction
assumptions, we have
$
\overline {s_{i+j+1}(\kappa)}\subset s_{i+j}(\rho)\subset s_i(\sigma)
$. This finishes the induction.

Let $A=\{\sigma\in {\cal T}_i^{(M+1)}\mid \sigma\subset\tau\}$.
Then $\tau = \bigcup_{\sigma
\in A} \sigma$, and
$$
g_i(\tau)=g_i(\bigcup_{\sigma
\in A} \sigma)=\bigcup_{\sigma
\in A} g_i(\sigma)\subset \bigcup_{\sigma
\in A} s_i(\sigma)\subset
\bigcup_{\sigma
\in A} \overline{s_i(\sigma)}\subset
s_L(\tau).
$$
Here the
last inclusion follows from what we proved above and the
first one holds by iii).

{\it Claim 1: For every $y\in B^{M+1}$ the limit 
${\rm lim}_{i\rightarrow \infty}g_i(y)$ exists.} 
Take an arbitrary open (in $X\cup \partial X$) finite cover 
${\cal W}$ of $\overline {U}\cap
\partial X$. Since for every  $j\ge i>0$
every $A\in U_j^1$ satisfies
$A\subset U\setminus B(x_0,N_{i-1})$, there exists a natural
$L>N_0$ such that for every 
$i\geq L$
every neighbourhood $A\in U_i^1$
is contained in 
some member of ${\cal W}$.
Take an arbitrary $y\in B^{M+1}$. Let $\tau \in {\cal T}_L$
be a maximal simplex containing $y$. 
Then, by Claim 0, $g_i(\tau)\subset s_L(\tau)\subset
W$ for every $i\geq L$ and some $W\in {\cal W}$. This implies
the existence of the limit.

{\it Claim 2: ${\rm lim}_{i\rightarrow \infty}
g_i(y)\in U\cap \partial X$.} Follows from: $g_i(y)
\in \bigcup U_1^1$ for every $i$; 
$\overline A\subset U$ for every $A\in U_1^1$.

{\it Claim 3: The formula
$g(x)={\rm lim}_{i\rightarrow \infty}g_i(x)$ 
defines a continuous map $g\colon B^{M+1}\rightarrow U\cap \partial X$.}
As in the proof of Claim 1, 
for every finite open cover ${\cal W}$ there exists $L>0$ such that
for every $i\geq L$ and any $A\in U_i^1$ the star
$\bigcup {\rm St}(A)$ of $A$ in $U_i^1$ is contained in some 
member of ${\cal W}$. 
Take an arbitrary $y\in B^{M+1}$. Let $\tau \in {\cal T}_L$
be a maximal simplex containing $y$. 
As in Claim 1, we have $g_i(\sigma)\subset s_L(\sigma)$ for every 
$i\geq L$ and every $(M+1)$-simplex $\sigma$ of ${\cal T}_L$
intersecting 
$\tau$,
and hence $g_i(\bigcup St(\tau))\subset \bigcup 
\lk s_L(\sigma) | \sigma \in St(\tau)\rk \subset \bigcup St(s_L(\tau))
\subset W$ for some $W\in {\cal W}$.
In other words, for every open cover ${\cal W}$ as above
and any given $y\in B^{M+1}$
there exists a natural $L$, $W\in {\cal W}$, and an open neighbourhood
$E \subset B^{M+1}$ of $y$ such that for every
$i\geq L$ we have $g_i(E)\subset W$. This implies that the limit
of $g_i$'s is continuous.

{\it Claim 4: The map $g\colon B^{M+1}\rightarrow U\cap \partial X$
extends $f$.} Follows from the fact that $g_i|_{S^M}=
p_{N_i}\circ h_i|{S^M}=p_{N_i}\circ f$
and ${\rm lim}_{i\rightarrow \infty}p_{N_i}\circ f(y)=f(y)$ for
every $y\in S^{M}$.
\hfill{$\diamond$}
\smallskip

\bf Proposition 3.5\sl\par
$\partial X$ is $(n-2)$-connected and locally $(n-2)$-connected.
\rm\par
Proof.
For the local statement let $x\in \partial X$ and let $W\ni x$
be its open (in $X\cup \partial X$) neighbourhood. By Lemma 2.1
one can find standard neighbourhoods $U$ and $V$ of $x$
contained in $W$ and such that $\overline V\subset U$. Then 
by Lemma 3.4 for
every $k\in \lk 0,1,...,n-2 \rk$ every map $f\colon S^k
=\partial B^{k+1}\rightarrow V$ 
has an extension $g\colon B^{k+1}\rightarrow U\subset W$.
For the global case apply Lemma 3.4 setting $V=U=X\cup \partial X$.
\hfill{$\diamond$}

\jump

\def\q{{\bf q}}

\def\ra{right-angled }
\def\st{star-like }
\def\Z{{\bf Z}}

\centerline{\bf 4. Right-angled buildings}
\medskip
Throughout this section $(W,S)$ is a finitely generated \ra Coxeter system,
not necessarily hyperbolic.
In subsection 4.1. we assume it to be finite,
i.e. $W\simeq(\Z/2)^n$, $S=\{(1,0,\ldots,0),\ldots,(0,\ldots,0,1)\}$.
\medskip

\bf 4.1. Finite \ra buildings.\rm
\smallskip
We will analyse the structure of finite $W$-buildings,
as well as maps between such buildings. 
This will be needed later for the constructions of infinite \ra buildings
and of maps between them. A typical step of those constructions 
consists of extending a map defined on a subset of a finite residue to
the whole residue.

We will treat a building
combinatorially, as a set (of chambers) equipped with  
a family $(\sim_s)_{s\in S}$ of equivalence relations 
(the adjacency relations). The standard example of 
a finite $W$-building is a \it product building\/\rm: 
the set of chambers $Y$ is a product
$\prod_{s\in S}Y_s$, where each $Y_s$ is a finite set of cardinality at least 2
(at least 3 if one wants a thick building). Two chambers
$(y_s), (y'_s)$ are $t$-adjacent if $y_s=y'_s$ for all $s\ne t$.
Apartments are of the form $A=\prod_{s\in S}A_s$, where 
each $A_s$ is a two-element subset of $Y_s$. 

It is easy to see that any $(\Z/2)^2$-building $X$ is a product building
(we will frequently apply this fact to residues in larger buildings). 
Indeed, let $S=\{s,t\}$ and let $Y_s=X/{\sim_s}$, $Y_t=X/{\sim_t}$.
Since  any two chambers $x,x'\in X$ are contained in some apartment,
$[x]_{\sim_s}$ and $[x']_{\sim_t}$ always have a common chamber.
Therefore, the map $X\ni x\mapsto([x]_{\sim_s},[x]_{\sim_t})\in Y_s\times Y_t$ is onto.
As no two chambers can be simultaneously $s$- and $t$-adjacent, this map is also injective.

By a morphism between two $W$-buildings (or subsets of such buildings) we mean a map of the sets
of chambers preserving the relations. A subset $E$ of a 
$W$-building $X$ is called {\it star-like} (with respect to a chamber $B\in X$)
if for every $x\in E$ every minimal gallery from $B$ to $x$ is contained in
$E$.
\smallskip
\bf Lemma 4.1\sl\par
Let $X$ be any finite $W$-building, and let $Y$ be a product $W$-building described above.
Let $E\subseteq X$ be \st with respect to a chamber $B$, and let $\psi\colon E\to Y$ be a morphism. 
Then $\psi$ extends to a morphism $\phi\colon X\to Y$. Moreover 
\item{1.} If two such extensions coincide on each class $[B]_{\sim_s}$, then they are equal.
\item{2.} If $\phi$ is injective on each class $[B]_{\sim_s}$, then $\phi$ is a monomorphism.
\item{3.} If $\phi$ maps each class $[B]_{\sim_s}$ onto $[\phi(B)]_{\sim_s}$, then $\phi$ is an epimorphism.
\rm\par
Proof.
Let $\pi\colon X\to W$ be the $B$-based folding map. Put $X_k=\pi^{-1}(\{w\in W\mid \ell(w)\le k\})$.
Let $C=\psi(B)$ if $B\in E$, or let $C$ be an arbitrary chamber of $Y$ if $E=\emptyset$.
Put $\phi(B)=C$. Define $\phi$ on $[B]_{\sim_s}\setminus (E\cup\{B\})$ to be an arbitrary map to
$[C]_{\sim_s}$; define $\phi$ on $[B]_{\sim_s}\cap E$ to coincide with the restriction of $\psi$;  
do this for each $s$. At this moment we have defined $\phi$ on $\bigcup_{s\in S}[B]_{\sim_s}=X_1$
so that it coincides with $\psi$ on $X_1\cap E$. Inductively on $k$ we will extend $\phi$ to $X_k$,
and check that the extension coincides with $\psi$ on $X_k\cap E$. Suppose this is done for $k-1\ge1$.
Let $x\in X$, $\pi(x)=w$, $\ell(w)=k$. For $t\in In(w)$ we denote by $x^t$ the chamber in the $t$-residue 
of $x$ which is closest to $B$. Since $x\sim_tx^t$, $\phi(x)$ has to be $t$-adjacent to $\phi(x^t)$ for every 
$t\in In(w)$.

Let $\phi(x^t)=(y^t_s)_{s\in S}$.
Let $t,t'\in In(w)$ be distinct, and let $s\in S$, $s\ne t,t'$. We claim that $y^t_s=y^{t'}_s$.
Indeed, let $x^{t,t'}$ be the chamber in the $\{t,t'\}$-residue of $x$ which is closest to $B$.
Since $W$ is  right-angled, we have $x^{t,t'}=(x^t)^{t'}=(x^{t'})^t$; consequently $x^{t,t'}\sim_tx^{t'}$, 
$x^{t,t'}\sim_{t'}x^t$. Therefore $y^t_s=\phi(x^{t,t'})_s=y^{t'}_s$.
Let us denote by $y_s$ the common value of  $y^t_s$,  $t\ne s$. Clearly, $y=(y_s)_{s\in S}$ is the unique chamber 
in $Y$ which is $t$-adjacent to $\phi(x^t)$ for each $t\in In(w)$. Therefore, we have to put $\phi(x)=y$.
Notice that if $x\in E$, then $x^t\in E$ for all $t\in In(w)$. Therefore $\psi(x)$
is $t$-adjacent to $\psi(x^t)=\phi(x^t)$ for all $t\in In(w)$, hence $\psi(x)=y=\phi(x)$. 
We apply the above procedure to every $x\in X_k\setminus X_{k-1}$, and get the required extension.

1. Follows from the construction: after defining $\phi$ on $X_1$ we made no choices.

2. Let $\pi_C\colon Y\to W$ be the $C$-based folding map. This map is given by
$\pi_C(y)=\Pi_{\{s\in S\mid C_s\ne y_s\}}s$. We first show that (under the assumption of 2)
the map $\phi$ is $\pi-\pi_C$ equivariant. Again, this is done by induction on $k$.
We have $\pi_C(\phi(B))=\pi_C(C)=1=\pi(B)$.
Then $\phi([B]_{\sim_s}\setminus \{B\})\subseteq[C]_{\sim_s}\setminus \{C\}$, $\pi([B]_{\sim_s}\setminus \{B\})=\{s\}=
\pi_C([C]_{\sim_s}\setminus \{C\})$, which checks $\pi$-equivariance on $X_1$.
Let now $x\in X$, $y=\phi(x)$, $w=\pi(x)$, $\ell(w)=k\ge2$. For $t\in In(w)$ we have $\pi(x^t)=wt$,
and, by the inductive assumption, $\pi_C(\phi(x^t))=wt$.
It follows that $\{s\in S\mid y^t_s\ne C_s\}=In(w)\setminus \{t\}$.
Hence $\{s\in S\mid y_s\ne C_s\}=In(w)$ and $\pi_C(\phi(x))=w$.

Consequently, if we have $\phi(x)=\phi(z)$, then $\pi(x)=\pi(z)$.
Let $x,z$ be such a pair with the shortest possible $w=\pi(x)$, and let $t\in In(w)$.
Then $\phi(x^t)=\phi(z^t)$, since both are the chamber 
in the $t$-residue of $\phi(x)=\phi(z)$ which is closest to $C$.
Since our counterexample to injectivity has shortest $w$, we deduce 
$x^t=z^t$---for all $t\in In(w)$. Notice that $\ell(w)\ge2$ 
($\pi$-equivariance and the assumption of 2 imply that $\phi$ is injective on $X_1$).
Let $t,t'$ be two distinct elements of $In(w)$. Then
$x\sim_tx^t=z^t\sim_tz$ and $x\sim_{t'}x^{t'}=z^{t'}\sim_{t'}z$, 
so that $x$ is both $t$- and $t'$-adjacent to $z$. This is possible in a building 
only if $x=z$.

3. By induction on $k$ we will prove that for any $(x,u)\in X_k\times S$ 
the map $\phi\colon [x]_{\sim_u}\to[\phi(x)]_{\sim_u}$ is surjective. The statement is true 
for $k=0$ by the assumption of 3.
Let $x\in X_k$, $\pi(x)=w$, $\ell(w)=k$, and let $u\in S$. We can assume that
$\ell(wu)=k+1$---otherwise $x^u\in X_{k-1}$, $[x]_{\sim_u}=[x^u]_{\sim_u}$ and the statement
for $(x,u)$ is true by the inductive assumption applied to $(x^u,u)$.
Pick a $t\in In(w)$. Let $y=(y_s)\in [\phi(x)]_{\sim_u}$. We then have $\phi(x)_s=y_s$ for $s\ne u$,
and $\phi(x^t)_s=y_s$ for $s\ne u,t$. Let $z_s=y_s$ for $s\ne t$, $z_t=\phi(x^t)_t$. Then
$z=(z_s)\sim_t \phi(x^t)$, hence (by the inductive assumption for $(x^t,u)$) there exists $x'\in[x^t]_{\sim_u}$ such that
$\phi(x')=z$. Observe that in the $\{u,t\}$-residue of $x^t$ there is a unique element
$x''$ which is $u$-adjacent to $x$ and $t$-adjacent to $x'$, while $y$ is the unique chamber 
in $Y$ which is $u$-adjacent to $\phi(x)$ and $t$-adjacent to $z$. Hence, $\phi(x'')=y$, where
$x''\in[x]_{\sim_u}$. 

It follows that the image of $\phi$ is closed under all adjacency relations, hence it is 
equal to $Y$.
\hfill{$\diamond$}

\smallskip
\bf Corollaries\par\rm
\item{1.} One can take $E=\bigcup_{s\in S}[B]_{\sim_s}$, choose
an arbitrary chamber $\psi(B)\in Y$ and for each $x\in [B]_{\sim_s}$ pick an arbitrary
$\psi(x)\in[\psi(B)]_{\sim_s}$; every such $\psi$ extends to a unique morphism.
\item{2.} Let $E=\bigcup_{s\in S}[B]_{\sim_s}$, let $Y_s=[B]_{\sim_s}$. 
Put $\psi(B)=(B)_{s\in S}$. For $x\sim_sB$ put $\psi(x)_t=B$ for $t\ne s$, 
$\psi(x)_s=x$. Then the extension $\phi\colon X\to Y$ is an isomorphism. 
Thus, any finite $W$-building is isomorphic to a product building.
Therefore, Lemma 4.1 holds with $Y$ replaced by any finite $W$-building.
\item{3.} A corollary of the proof: every monomorphism of finite $W$-buildings $\phi\colon X\to Y$ is
$\pi$-equivariant (where $\pi\colon X\to W$ is a folding map based at an arbitrary chamber 
$x\in X$, and $\pi\colon Y\to W$ the $\phi(x)$-based folding map).
\medskip

\bf 4.2. Maps of infinite \ra buildings.\rm
\smallskip
\bf Definition\rm\par
A \it standard $W$-building \rm is a set $X$ (of chambers) equipped with: (a) a 
family $(\sim_s)_{s\in S}$ of equivalence relations with finite 
equivalence classes; (b) a morphism $\pi\colon X\to W$, called the folding map; 
such that the following are satisfied:

\item{1.} $(\forall x\in X)(\forall s\in S)(\exists x'\in X)(x\sim_sx'\land x\ne x')$;
\item{2.} $\pi^{-1}(1)$ has one element (denoted $B$ and called the base chamber);
\item{3.} let $x\in X$, $T=In(\pi(x))$, $w=\pi(x)$.
Then $Res(x,T)$ is a finite \ra building and the map $Res(x,T)\ni x'\mapsto (ww_T)^{-1}\pi(x')\in W_T$ is
a folding map of that building (where $w_T$ is the longest element in $W_T$).
\smallskip
It is pretty clear that any locally finite $W$-building with any folding map
is a standard $W$-building. In particular, condition 3 follows from property (F3) 
stated in section 1. More specifically, $Res(x,T)$ is mapped by $\pi$ onto 
the coset $wW_T$ of $W_T$; $w=\pi(x)$ is the longest element of $wW_T$,
therefore $ww_T$ is the shortest element of $wW_T$. Let $y$ be the shortest chamber
in $Res(x,T)$ (as in (F3)); then $\pi(y)=ww_T$. The $y$-based folding map of $Res(x,T)$
is the composition of (restricted) $\pi$ and the left multiplication in $W$ that
moves the coset $wW_T$ to $W_T$ and the element $\pi(y)$ to $1$. 
This left multiplication is the left multiplication by $\pi(y)^{-1}=(ww_T)^{-1}$.
\smallskip

\bf Remarks\rm\par
\item{1.} Later we will prove that a standard building is in fact a building. 
\item{2.} The residue $Res(x,T)$ in condition 3 intersects $\pi^{-1}(ww_T)$ in one chamber, to be called 
the shortest chamber of $Res(x,T)$. The folding map in condition 3 is based at that chamber.  
\item{3.} Conditions 1 and 3 together imply that for every $x\in X$ and every spherical
$T\subseteq S$ the residue $Res(x,T)$ is a finite $W_T$-building, and the restriction of
$\pi$ composed with left multiplication by the inverse of the shortest element 
of $\pi(Res(x,T))$ is a folding map of that building. 
\item{4.} It follows from the previous remark that if $t\in In(\pi(x))$ then
$\pi^{-1}(\pi(x)t)\cap [x]_{\sim_t}$ consists of a unique element (to be denoted $x^t$). 
\smallskip
\bf Definition\rm\par
A \it local $W$-building \rm is a set $Y$ (of chambers) equipped with  a 
family $(\sim_s)_{s\in S}$ of equivalence relations, such that: 
\item{1.} for every $y\in Y$ and every spherical $T\subseteq S$, 
$Res(y,T)$ is a finite $W_T$-building;
\item{2.} $Y$ is \it gallery connected, \rm i.e.,
for every $y,y'\in Y$ there exists a gallery from $y$ to $y'$: a sequence $y_0=y, y_1,\ldots, y_k,y_{k+1}=y'$, such that
$y_i\sim_{s_i}y_{i+1}$ for some $s_i\in S$, where $i=0,1,\ldots,k$. 

\smallskip
The following theorem is rather weak. The proof will give us an idea of what should really be done.
\smallskip
 
\bf Theorem 4.2\sl\par
Let $X$ be a standard $W$-building, and let $Y$ be a local $W$-building. 
Then there exists a morphism $\phi\colon X\to Y$.\rm
\smallskip

Proof. 
Choose a well-ordering $<$ on $X$ such that each initial segment $X_{<x}$ is \st (with respect to
$B$). We define $\phi$ inductively. To start, we pick  any $y\in Y$ and declare $\phi(B)=y$.
At limit steps we take union. Suppose $\phi\colon X_{<x}\to Y$ has already been defined.
Let $T=In(\pi(x))$. Since $X_{<x}$ is star-like, so is $X_{<x}\cap Res(x,T)$ (in $Res(x,T)$, with respect
to the shortest element $x_0$ of that residue). Since $\phi$ is a morphism,
it maps $X_{<x}\cap Res(x_0,T)$ into $Res(\phi(x_0),T)$; this restriction can, by Lemma 4.1, be extended
to $\eta\colon Res(x_0,T)\to Res(\phi(x_0),T)$. We put $\phi(x)=\eta(x)$. Since all chambers in $X_{<x}$ which
are adjacent to $x$ belong to $Res(x_0,T)$ (where $\phi$ coincides with $\eta$), 
the extended $\phi\colon X_{<x}\cup\{x\}\to Y$ is a morphism.
\hfill{$\diamond$}
\smallskip
Notice that, in the construction of $\phi$ described in the proof, if $In(\pi(x))$
has at least two elements, then $\phi(x)$ is uniquely determined by $\phi|_{X_{<x}}$.
In fact, if $u,t\in In(\pi(x))$, then $\phi(x)$ is uniquely determined by 
$\phi(x^t)$ and $\phi(x^u)$: it is 
the unique chamber $u$-adjacent
to $\phi(x^t)$ and $t$-adjacent to $\phi(x^u)$ (as in the proof of Lemma 4.1).
If, on the other hand, $In(\pi(x))=\{s\}$, then $\phi(x)$
can be freely chosen in $[\phi(x_0)]_{\sim_s}$. These observations are basic for the next theorem.

Let $X$ be a standard $W$-building, and let $Y$ be a local $W$-building.  
We say that a morphism $\phi\colon X\to Y$ is a \it local monomorphism (resp. local epimorphism, covering map),
\rm if for every $x\in X$ and every spherical $T\subseteq S$ the residue $Res(x,T)$ is injectively 
(resp. surjectively, bijectively) mapped by $\phi$ to $Res(\phi(x),T)$. 
\smallskip
\bf Definition\rm\par
The \it root set \rm $R(X)$ of a standard $W$-building $X$ is 
$\{(x,s)\in X\times S\mid In(\pi(x)s)=\{s\}\}$.

\smallskip

\bf Theorem 4.3\sl\par
Let $X$ be a standard $W$-building, let $Y$ be a local $W$-building, 
and let $\phi\colon X\to Y$ be a morphism. 
Let $R$ be the root set of $X$.
\item{1.} The map $\phi$ is uniquely determined by $\phi(B)$ and the restrictions 
of $\phi$ to $[x]_{\sim_s}$, over all $(x,s)\in R$.
\item{2.} If all the above restrictions are injective, then $\phi$ is a local monomorphism.
\item{3.} If, for each $(x,s)\in R$, $\phi$ maps $[x]_{\sim_s}$ onto
$[\phi(x)]_{\sim_s}$, then $\phi$ is a local epimorphism and a surjection.
\item{4.}  If, for each $(x,s)\in R$, $\phi$ maps $[x]_{\sim_s}$ bijectively onto
$[\phi(x)]_{\sim_s}$, then $\phi$ is a covering map.\rm
\smallskip

Proof.

1. Let $\phi_1,\phi_2\colon X\to Y$ coincide on $B$ and on each $[x]_{\sim_s}$, $(x,s)\in R$.
Suppose that $x\in X$ is a chamber with shortest $w=\pi(x)$ such that $\phi_1(x)\ne\phi_2(x)$.
If $In(w)=\emptyset$, then $w=1$ and $x=B$---contradiction.
If $In(w)=\{s\}$, then $(x^s,s)\in R$, $x\in[x^s]_{\sim_s}$---contradiction again.
If $T=In(w)$ has at least two elements, then, by Lemma 4.1 part 1 applied to $Res(x,T)$, 
$\phi_1(x)$ and $\phi_2(x)$ are uniquely determined by $\phi_1|_{X_k}=\phi_2|_{X_k}$ (where $k=\ell(w)-1$),
so that they coincide. 

2. Suppose not. Let $x_0\in X$ be an element of $X$ with the shortest possible $w=\pi(x_0)$, such that 
for some spherical $T$ the restriction of $\phi$ to $Res(x_0,T)$ is not injective.
By Lemma 4.1 part 2, there exists $t\in T$ and chambers $x,x'\sim_tx_0$ such that $\phi(x)=\phi(x')$.
If $In(wt)\ne \{t\}$, then  $Res(x,In(wt))$ is a residue on which 
$\phi$ is injective, and whose shortest chamber is shorter than $x_0$, contradiction.
In the case $In(wt)=\{t\}$ we have $(x_0,t)\in R$ and $x,x'\in[x_0]_{\sim_t}$, so that $\phi(x)\ne \phi(x')$, contradiction.   
   
3. Suppose that $\phi$ is not a local epimorphism. Let $Res(x_0,T)$ be a counterexample with shortest
$w=\pi(x_0)$. Then, by Lemma 4.1 part 3, there is a $t\in T$ such that $\phi\colon [x_0]_{\sim_t}\to[\phi(x_0)]_{\sim_t}$
is not onto. As in the proof of 2 we see that $In(wt)=\{t\}$. Therefore $(x_0,t)\in R$, contradiction.

Since the image of a local epimorphism is closed under the adjacency relations, and since $Y$ is gallery connected,
we have $\phi(X)=Y$.

4. Follows from 2 and 3.
\hfill{$\diamond$}
\medskip

\bf 4.3. Construction.\rm
\smallskip
In this subsection we present a construction of a general standard $W$-building. 
 Let $W=\{w_1=1,w_2,\ldots\}$ be a numbering of elements of $W$ such that each $W_k=\{w_1,\ldots,w_k\}$ is 
a \st subset of $W$ (with respect to $1$). The process of building $X$ is inductive.
At $k$'th step we construct the part $X_k$ of $X$ which is going to be the preimage of $W_k$
under the folding map. To get $X_k$ from $X_{k-1}$ we need to attach the chambers that fold 
to $w_k$. Such a chamber $x$ is contained in a finite residue $Res(x,In(w_k))$ which is isomorphic
to a product building and whose large part is contained in $X_{k-1}$.
Thus, $X_k$ is obtained from $X_{k-1}$ by gluing to it product buildings that will become
$Res(x,In(w_k))$ for $x\in \pi^{-1}(w_k)$.  

We now proceed to the details.
We would like to construct, by induction on $k$, sets $X_k$ with equivalence relations
$(\sim^k_s)_{s\in S}$, together with morphisms $\pi_k\colon X_k\to W_k$, such that:
\item{1.} $X_{k-1}\subseteq X_k$. 
\item{2.} Each relation $\sim^k_s$ when restricted from 
$X_k\times X_k$  to
$X_{k-1}\times X_{k-1}$ yields the relation $\sim^{k-1}_s$. (Therefore, we will simply use $\sim_s$). 
\item{3.} $\pi_k|_{X_{k-1}}=\pi_{k-1}$. (Again, we will often denote the map $\pi_k$ simply by $\pi$). 
\item{4.} $\pi_k^{-1}(W_{k-1})=X_{k-1}$.
\item{5.} Let $x\in X_k$, $T=In(\pi(x))$, $w=\pi(x)$.
Then $Res_k(x,T)$ is a finite \ra building and the map $Res_k(x,T)\ni x'\mapsto (ww_T)^{-1}\pi(x')\in W_T$ is
a folding map of that building. Here $Res_k$ stands for the residue in $X_k$. 

\noindent Finally, we will obtain a standard $W$-building $X=\bigcup_k X_k$ with the folding map $\pi=\bigcup_k \pi_k$. 
In fact, at $k$th step we will construct not only $X_k$ and $\pi_k$, but also the following additional data:
\item{(a)} an integer $q_{x,s}\ge1$, for each $(x,s)\in X_k\times S$;

\item{(b)} for each $u\in W$ such that $uw_U\in W_k$ (where $U=In(u)$) and each $y\in\pi^{-1}(uw_U)$:
a $(uw_U)^{-1}\pi_k$-$\pi_{y,U}$ equivariant monomorphism  $\phi_{y,U}\colon  Res_k(y,U)\to Y_{y,U}$. 
Here 
$Y_{y,U}$ is a product $W_U$-building with the $s$-factor $Y_{y,s}$ of cardinality
$q_{y,s}+1$, and $\pi_{y,U}\colon Y_{y,U}\to W_U$ is the $\phi_{y,U}(y)$-based folding map.
We will usually briefly say that $\phi_{y,U}$ is $\pi$-equivariant.

\noindent The numbers $q_{x,s}$ are subject to  extra conditions:
\item{6.} If $q_{z,s}$ and $q_{z',s}$ are defined and $z\in Res_k(z',T)$ for a spherical $T\subseteq S$ containing $s$,
then $q_{z,s}=q_{z',s}$.
\item{7.} If $y\in X_k$ and $\{\pi(y),\pi(y)s\}\subseteq W_k$, then
$Res_k(y,s)$ has $q_{y,s}+1$ elements.  
\smallskip

The first step is: $X_1=\{B\}$. We choose the numbers $q_{B,s}$ and the maps $\phi_{B,T}\colon \{B\}\to Y_{B,T}$
arbitrarily.

Suppose that we have already constructed everything promised for $k-1$. Let $w=w_k$, let
$T=In(w)$. The set $X_k$ is obtained from $X_{k-1}$ by gluing
$Y_{x,T}$, for all $x\in\pi^{-1}(ww_T)$, via the maps $\phi_{x,T}$. 
Throughout the proof, $x$ will be a generic notation  for an element of $\pi^{-1}(ww_T)$.

\bf Lemma\sl

Let $x_1,x_2\in \pi^{-1}(ww_T)$, $x_1\ne x_2$. Then $Res_{k-1}(x_1,T)\cap Res_{k-1}(x_2,T)=\emptyset$.
\par\rm
Proof. Suppose not; then $x_2\in Res_{k-1}(x_1,T)$. The map $\phi_{x_1,T}$ being $\pi$-equivariant, we have
$\phi_{x_1,T}(x_2)=\phi_{x_1,T}(x_1)$. However, $\phi_{x_1,T}$ is injective; hence 
$x_1=x_2$, contradiction.
\null\hfill{$\diamond$(Lemma)}

\item{1.} Follows from the fact that the gluing maps are injective.

\item{2.} \it Claim\/\rm: If $y,z\in Res_{k-1}(x,T)$ and $\phi_{x,T}(y)\sim_t\phi_{x,T}(z)$, then $y\sim_tz$.
\item{}
\it Proof. \rm We can assume that $t\in T$, for otherwise $y=z$.
Let $\phi=\phi_{x,T}$, $Y=Y_{x,T}$. There are two cases.

\itemitem{a)} {$\{\pi(y),\pi(y)t\}\subseteq W_{k-1}$. In that case, due to 7., $Res_{k-1}(y,t)$ has $q_{y,t}+1$
elements; $Res_Y(\phi(y),t)$ has $q_{x,t}+1$ elements. However, 6. implies that $q_{x,t}=q_{y,t}$ ($y\in Res_{k-1}(x,T)$); 
therefore $\phi$ restricts to a bijection between these residues.
Since $\phi$ is an injection, this implies that $z\in Res_{k-1}(y,t)$.}
\item{}
\vskip-.4cm
\itemitem{b)} {$\pi(y)t=w$. Then $Res_Y(\phi(y),t)$ consists of $\phi(y)$ and chambers which map to $w_T$ under $\pi_{x,T}$.
Hence either $\phi(z)=\phi(y)$ (and---$\phi$ being injective---we get $z=y$) or
$\pi_{x,T}(\phi(z))=w_T$. In the latter case $\pi$-equivariance of $\phi$ implies
$\pi(z)=w$, which contradicts $z\in X_{k-1}$.\hfill{$\diamond$({\it Claim})}
\item{}
Due to the claim, the following definition makes sense:
$y\sim^k_sz$ if either $y,z\in X_{k-1}$ and $y\sim^{k-1}_sz$, or 
$y,z\in Y_{x,T}$ for some $x$ and $y\sim_sz$ in $Y_{x,T}$.
\item{}
Finally, we need to check that $\sim_s^k$ is an equivalence relation, 
the only non-trivial condition being transitivity:
$(a\sim_sb\land b\sim_s c)\Rightarrow a\sim_sc$. The cases $a,b,c\in X_{k-1}$ and
$a,b,c\in Y_{x,T}$ are clear. Thus, we can assume that at least one of $a,b,c$ is  in 
$Y_{x,T}\setminus X_{k-1}$ 
(for some $x$). Then we can assume that $s\in T$ (otherwise $a=b=c$). 
If $b\in Y_{x,T}\setminus X_{k-1}$, then $a,c\in Y_{x,T}$ and $a\sim_sc$ follows.
If not, we can assume $a\in Y_{x,T}\setminus X_{k-1}$, $b\in Y_{x,T}\cap X_{k-1}=Res_{k-1}(x,T)$.    
Now if $c\in X_{k-1}$, then $c\sim_sb$ implies $c\in Res_{k-1}(x,T)\subseteq Y_{x,T}$,
and $a\sim_sc$ follows. If $c\not\in X_{k-1}$, then $c\in Y_{x',T}$ (for some $x'\in\pi^{-1}(w)$).
It follows that $b\in Y_{x',T}\cap Y_{x,T}$, hence, in view of the Lemma, $x=x'$ and $a,b,c\in Y_{x,T}$.
 
\item{3.} 
We define $\pi_k$ as follows: if $y\in X_{k-1}$, then $\pi_k(y)=\pi_{k-1}(y)$; 
if $y\in Y_{x,T}$ we put $\pi_k(y)=ww_T\pi_{x,T}(y)$. This definition is correct because of
the $\pi$-equivariance of $\phi_{x,T}$. Condition 3. is clear.

\item{4.}
Suppose that $y\in Y=Y_{x,T}$, but $y$ is not in the image of $\phi=\phi_{x,T}$.
We claim that $\pi(y)=w$, or, equivalently, that $\pi_{x,T}(y)=w_T$.
Suppose not; let $y$ be a counterexample with shortest $u=\pi_{x,T}(y)$.  
Notice that $u\ne 1$, because $\pi_{x,T}^{-1}(1)=\{\phi(x)\}$.
Let $t\in In(u)$, and let $y^t=\phi(z)$. As in 2., we have
$q_{x,t}=q_{z,t}$. Moreover, $\pi(z)=ww_Tut$ and $\pi(z)t=ww_Tu$ belong to $W_{k-1}$,
so that $Res_{k-1}(z,t)$ has cardinality $q_{z,t}+1$---the same as $Res_Y(y^t,t)$. Therefore $\phi$ maps $Res_{k-1}(z,t)$
bijectively onto $Res_Y(y^t,t)$, and $y$ is in the image of $\phi$, contradiction.

\item{5.}
The new residues to be checked are $Res_k(y,T)$, for $y\in Y_{x,T}$, $\pi(y)=w$.
But in this case $Res_k(y,T)=Y_{x,T}$, and $(ww_T)^{-1}\pi=\pi_{x,T}$.

\item{6.}
{Let $y\in Y_{x,T}\setminus X_{k-1}$, $s\in S$. If $s\in T$, we put $q_{y,s}=q_{y^s,s}$. 
If $s\not\in T$, but there exists $t\in T$ such that
$\{t,s\}$ is spherical, then we put $q_{y,s}=q_{y^t,s}$. This does not lead to contradictions:
if $t'\in T$ and $\{s,t'\}$ is spherical, then $\{s,t,t'\}$ is also spherical, 
and $y^{t'}\in Res_{k-1}(y^t,\{s,t,t'\})$, so that $q_{y^{t'},s}=q_{y^t,s}$ by 6.
Finally, if $s\not\in T$ and no $t\in T$ commutes with $s$, then we choose $q_{y,s}$ arbitrarily.
Observe that this last case occurs exactly when $In(ws)=\{s\}$.

Now suppose that $z,z'\in X_k$, $z\in Res_k(z',U)$ where $U$ is spherical and $s\in U$.
We will show that $q_{z,s}=q_{z',s}$. Let $z=z_1,z_2,\ldots,z_m=z'$ be a gallery in $X_k$, 
$z_i\sim_{u_i}z_{i+1}$, $u_i\in U$. Suppose that some two consecutive chambers $z_i,z_{i+1}$ 
do not belong to $X_{k-1}$. Then they are both in $Y_{x,T}$ for some $x$, and $u_i\in T$.
We insert $z^{u_i}_i$ between $z_i$ and $z_{i+1}$. Repeating the process we ensure that 
if $z_i\not\in X_{k-1}$, then $z_{i-1},z_{i+1}\in X_{k-1}$ (except $i=0,m$).
Then we replace each triple $z_{i-1},z_i,z_{i+1}$ with $z_i\not\in X_{k-1}$ 
by $z_{i-1},z^{u_i}_{i-1}=z^{u_{i-1}}_{i+1},z_{i+1}$. We obtain a $U$-gallery from $z$ to $z'$
whose all but external chambers lie in $X_{k-1}$. We conclude that $q_{z,s}=q_{z_1,s}=q_{z_{m-1},s}=q_{z',s}$.}

\item{7.} The new residues to be checked are $Res_k(y,s)$ for $y$ such that $\pi(y)=w$ or $\pi(y)=ws$
(where---for some $x$---$y\in Y_{x,T}$ and $s\in T$). In either case, $Res_k(y,s)=Res_Y(y,s)$ has $q_{x,s}+1$ elements.
However, $y\in Res_k(x,T)$ so that, by 6., $q_{x,s}=q_{y,s}$. 

\smallskip

Finally, we need to construct or extend some of the maps $\phi_{y,U}$.  

Some cases are easy. If $u=w$, then $U=T$ and 
$y\in \pi^{-1}(ww_T)$. Then $Res_k(y,T)=Y_{y,T}$ and we put $\phi_{y,T}=Id_{Y_{y,T}}$.
If $uw_U=w$, then we choose $\phi_{y,U}\colon \{y\}\to Y_{y,U}$ arbitrarily.
If $w\not\in uW_U$, then we do not change $\phi_{y,T}$.

Thus, we can assume that $w\in uW_U$, but $w\ne uw_U$.
Then $w=uw_Uu_1\ldots u_k$ for some pairwise different $u_1\ldots u_k\in U$.
We have $u_i\in In(w)=T$, so that $w\in uw_UW_{U\cap T}$, or equivalently,
$uw_U\in wW_{U\cap T}$. Since $uw_U$ is the shortest element in $uW_U=wW_U$,
it is also the shortest element in $wW_{U\cap T}$.

\it Claim\/: \rm There is an $x\in \pi^{-1}(ww_T)$ such that $y\in Y_{x,T}$.

\it Proof. \rm Since $uw_U\in wW_{U\cap T}\subseteq wW_T= ww_TW_T$, there is a $T'\subseteq T$ such that 
$uw_U=ww_Tw_{T'}$. Moreover, $T'\subseteq In(uw_U)$, therefore
$ww_T\in uw_UW_{T'}\subseteq uw_UW_{In(uw_U)}=\pi(y)W_{In(\pi(y))}$. By 5.,
there is an $x\in \pi^{-1}(ww_T)\cap Res_{k-1}(y,T')$ (the folding map, when restricted to a residue,
is onto a suitable coset of the corresponding special subgroup of the Coxeter group).
Then $y\in Res_k(x,T')\subseteq Res_k(x,T)=Y_{x,T}$.
\hfill{$\diamond$({\it Claim})}

Let $Y_{y,U\cap T}=Res_{Y_{x,T}}(y,U\cap T)$.
Since $Res_{k-1}(y,U)\cap Y_{y,U\cap T}=Y_{y,U\cap T}\setminus \pi^{-1}(w)$,
it is a \st set in the building $Y_{y,U\cap T}$ (with respect to $y$), therefore $\phi_{y,U}$
extends to a monomorphism $\psi\colon Y_{y,U\cap T}\to Y_{y,U}$.
Gluing $\phi_{y,U}$ with $\psi$ we get an extended map
$\phi_{y,U}\colon Res_{k-1}(y,U)\cup Y_{y,U\cap T}\to Y_{y,U}$.
We claim that this map is injective:
indeed, $\psi$ is $\pi$-equivariant as a monomorphism of buildings, and
hence the extended map is $\pi$-equivariant. Furthermore, if $z\in Res_{k-1}(y,U)$
and $z'\in Y_{y,U\cap T}\setminus Res_{k-1}(y,U)$, then $\pi(z)\ne w=\pi(z')$, therefore
$\phi_{y,U}(z)\ne\phi_{y,U}(z')$. Now it is enough to observe that $\phi_{y,U}$ is injective
on $Res_{k-1}(y,U)$ and that $\psi$ is injective.

Finally, we claim that $Res_{k-1}(y,U)\cup Y_{y,U\cap T}=Res_k(y,U)$.
Since $\pi$ is a morphism, we have that $Res_k(y,U)\subseteq \pi^{-1}(uw_UW_U\cap W_k)$.
Any $U$-gallery in $X_k$ starting at $y$ and ending at $z\not\in\pi^{-1}(w)$ can be 
modified, using the technique from the proof of 6., to a $U$-gallery not containing
chambers from $\pi^{-1}(w)$. This means that $Res_k(y,U)\cap X_{k-1}=Res_{k-1}(y,U)$.
Suppose now that $z\in \pi^{-1}(w)\cap Res_k(y,U)$. Then 
$Res_k(z,U\cap T)\subseteq Res_k(y,U)$,
$Res_k(z,U\cap T)\setminus \pi^{-1}(w)\subseteq Res_{k-1}(y,U)$. In particular, $Res_k(z,U\cap T)$
has a unique shortest element, lying in $\pi^{-1}(uw_U)\cap Res_{k-1}(y,U)$.
But the latter set equals $\{y\}$, because $\phi_{y,U}\colon Res_{k-1}(y,U)\to Y_{y,U}$
is a $\pi$-equivariant monomorphism. Therefore $y\in Res_k(z,U\cap T)$, $z\in Y_{y,U\cap T}$.

\smallskip
\bf Remark\rm

Notice that we were free to choose $q_{x,s}$ exactly for 
$(x,s)$ in the root set of $X$. 
\medskip

\bf 4.4. Uniqueness and lattices.\rm
\smallskip 

\bf Theorem 4.4\sl\par
For every local $W$-building $Y$ and any chamber $y\in Y$ 
there exists a standard $W$-building $X$ and a covering map $\phi\colon X\to Y$, $\phi(B)=y$.
\rm
\smallskip

Proof.
We perform the construction of $X$ as in subsection 4.3, together with the construction of 
$\phi$ as in the proof of Theorem 4.2. We put $X_1=\{B\}$ and $\phi(B)=y$.
Whenever we construct a chamber $x$ such that a choice of $q_{x,s}$ is needed for some $s$,
we choose $q_{x,s}=|[\phi(x)]_{\sim_s}|-1$. The pair $(x,s)$ will belong to the root set of $X$.
Later, when $\pi^{-1}(\pi(x)s)$ is constructed, we  are free to choose $\phi\colon [x]_{\sim_s}
\to[\phi(x)]_{\sim_s}$ (extending $x\mapsto\phi(x)$); we choose a bijection. By Theorem 4.3 part 4
we obtain a covering map.
\hfill{$\diamond$}
\smallskip

To talk about universal covers it is convenient to switch to the topological category
(and back). The geometric realisation $|X|$ of a local $W$-building $X$ 
is the geometric realisation of the poset of finite type residues in $X$ (i.e. $T$-residues for all spherical $T$).
One can label each vertex in $|X|$ with the type of the corresponding residue.
Let us give another description of $|X|$. 
Let $L$ be the finite simplicial complex with vertex set $L^{(0)}=S$,
a set of generators spanning a simplex in $L$ if and only if they pairwise commute.
We denote by $L'$ the first barycentric subdivision of $L$,
and by $CL'$ the cone over $L'$.
Then $|X|$ is $X\times CL'/\sim$, where
$(x,p)\sim(x',p')\iff p=p'$ and $x'\in Res(x,S(p))$;
here $S(p)=\{s\in L^{(0)}\mid \exists \sigma\in L', s\in\sigma,p\in|\sigma|\}$.
If $p\colon \widetilde{|X|}
\to|X|$ is any covering of $X$, then $\widetilde{|X|}$ is in fact a geometric realisation of 
a local $W$-building $\widetilde{X}$:  $\widetilde{X}$ is the preimage by $p$ of the
set of $\emptyset$-labelled vertices of $|X|$; $\widetilde{x}\sim_s\widetilde{x}'\iff
$ there exists a vertex $v\in\widetilde{|X|}$, joined by edges to $\widetilde{x}$
and to $\widetilde{x}'$, and such that $p(v)$ is of type $\{s\}$.
Notice that for any spherical $U\subseteq S$, and any residue $R$ in $X$ of type $U$,
the set $|X_{\le R}|$ is contractible (as a cone with apex $R$), hence its preimage by $p$
is a disjoint union of its homeomorphic copies. Therefore, $\widetilde{X}$ is indeed
a local $W$-building. We say that $\widetilde{X}$ is the universal cover of 
$X$ if $\widetilde{|X|}$ is the universal cover of $|X|$.
It is clear that morphisms of local $W$-buildings 
induce simplicial label-preserving
maps, and covering maps of local $W$-buildings induce simplicial covering maps.
\smallskip

\bf Theorem 4.5\sl\par
A standard $W$-building is a building.
\rm\smallskip
Proof.

\bf Lemma A\sl\par
The geometric realisation of a  standard $W$-building is contractible.
\rm\par
Proof. Follows Serre's proof for buildings [Se].
Let $W=\{w_1=1,w_2,\ldots\}$ be a numbering of elements of $W$ such that each set $W_k=\{w_1,\ldots,w_k\}$ is star-like.
Let $X_k=\pi^{-1}(W_k)$. The strategy is to show that $|X_k|$ deformation retracts onto $|X_{k-1}|$. 
To do this, it is enough to check  that each chamber $|x|$ in $|X_k\setminus X_{k-1}|$ deformation 
retracts onto $|x|\cap|X_{k-1}|$.
Let $T=In(w_k)$; then the pair $(|x|,|x|\cap|X_{k-1}|)$ is isomorphic to 
$(K,K^T)$, where $K^T=\bigcup_{t\in T}K_t$. 
Since $K$ is a cone over $K^S$, it is contractible.
It is checked in 
[D1] that $K^T$ is contractible for all spherical $T$. 
It follows that $K$ deformation retracts onto $K^T$.
\hfill{$\diamond$(Lemma A.)}
\smallskip
Let $X$ be a standard $W$-building.
\smallskip

\bf Lemma B\sl\par
Suppose that $\sigma\colon W\to X$ is a local monomorphism, such that 
$\sigma(1)=B$. Then $\sigma$ is a section of $\pi$
(i.e., $\pi\circ\sigma=Id_W$).
\rm\par
Proof.
We argue by contradiction. Let $w\in W$ be the shortest element such that $\pi(\sigma(w))\ne w$.
Let $s\in In(w)$; then $\pi(\sigma(ws))=ws$.
The map $(ws)^{-1}\pi\colon  Res(\sigma(ws),s)\to W_{\{s\}}$ is a folding map
(cf. remark 3 after the definition of a standard building), so that
$\pi(\sigma(ws))=ws$ while for any  $x\in Res(\sigma(ws),s)\setminus \{\sigma(ws)\}$
we have $\pi(x)=w$. Since $w\in Res(ws,s)$, we have $\sigma(w)\in Res(\sigma(ws),s)$;
but $\sigma$ is injective on $Res(ws,s)$ so that $\sigma(w)\ne \sigma(ws)$. 
Therefore $\pi(\sigma(w))=w$, contradiction.
\hfill{$\diamond$(Lemma B)}
\smallskip  

We define an apartment in $X$ as the image of any monomorphism
$\sigma\colon W\to X$. Notice that if $B\in\sigma(W)$, then we can modify 
$\sigma$ by precomposing it with left multiplication by $\sigma^{-1}(B)$, so
as to have $\sigma(1)=B$. Therefore, every apartment containing $B$ is the image of a 
section of $\pi$.
Observe that if $\sigma$, $\sigma'$ are sections of $\pi$ and $\sigma(W)$,
$\sigma'(W)$ are two apartments containing $B$ and a chamber $x$, then
$\sigma'\circ\pi\colon \sigma(W)\to \sigma'(W)$ is an isomorphism fixing $B$ and $x$.

Recall that $x^t$ denotes the shortest element in $Res(x,t)$ (we use this notation only if $t\in In(\pi(x))$).
Define inductively $x^{t_1\ldots t_it_{i+1}}=(x^{t_1\ldots t_i})^{t_{i+1}}$.
Notice that if $st=ts$ then $x^{st}=x^{ts}$, since both chambers are equal to the 
shortest element in $Res(x,\{s,t\})$. We will also use the fact that if $\sigma(w)=x$ for a section 
$\sigma$ of $\pi$, then $\sigma(wt)=x^t$ (assuming $t\in In(w)$).

\smallskip

\bf Lemma C\sl\par
For any $x\in X$ there exists a morphism $\sigma\colon W\to X$ which
is a section of $\pi$ and satisfies $\sigma(\pi(x))=x$.
\rm\par
Proof.
Induction on the length $k$ of $\pi(x)$. 

For $k=0$ we have $x=B$, and we just have to show the existence of a section of $\pi$.
A morphism $\sigma\colon W\to X$ such that $\sigma(1)=B$ can be constructed as in 
Theorem 4.2. Moreover, since $[x]_{\sim_s}$ has at least 2 elements for each 
$x\in X$, $s\in S$, there exists $\sigma$ which is injective on $[x]_{\sim_s}$ for each
$(x,s)\in R(W)$,  hence (by Theorem 4.3.2) is a local monomorphism. 
Then, by Lemma B, $\sigma$ is a section.

Now let $k>0$. Then we can find $w$ of length $k-1$ such that $\pi(x)=ws$ for some $s\in T=In(ws)$, 
and a section $\xi\colon W\to X$ of $\pi$ such that $\xi(w)=x^s$.
Let $us$ be the shortest element in $H(w,s)$, and let $ws,wst_1,\ldots,wst_1\ldots t_k=us$ be a minimal
gallery with $t_i\in \{s\}'$ (cf. Lemma 2.5).
We choose a well-ordering on $W$ with \st initial segments such that all elements 
of $H(w,s)$ are larger than all other elements. 
If $\{g,gt\}\cap H(w,s)=\emptyset$ and $(g,t)\in R(W)$ then we put $\sigma(gt)=\xi(gt)$
(so that $\sigma$ and $\xi$ coincide on $W\setminus H(w,s)$). Then we put $\sigma(us)
=x^{t_1\ldots t_k}$, and afterwards we only care about making injective choices, so
as to keep $\sigma$ a local monomorphism (and hence, as in the case $k=0$, a section of $\pi$).
We claim that $\sigma(ws)=x$. To check this we prove by descending induction on $i$  that 
$\sigma(wst_1\ldots t_i)=x^{t_1\ldots t_i}$. Indeed, observe that $\sigma(wt_1\ldots t_i)=
\xi(wt_1\ldots t_i)=x^{st_1\ldots t_i}=x^{t_1\ldots t_is}$, while by the inductive assumption
$\sigma(wst_1\ldots t_{i+1})=x^{t_1\ldots t_{i+1}}$. Since $x^{t_1\ldots t_i}$ is the unique 
chamber which is respectively $s$- and $t_{i+1}$-adjacent to the above two chambers,
it has to be equal to $\sigma(wst_1\ldots t_i)$.
\hfill{$\diamond$(Lemma C)}
\smallskip

For $x\in X$ we can find a standard $W$-building $X'$ and a covering
map $\phi\colon X'\to X$, $\phi(B')=x$ (Theorem 4.4).
By Lemma A, $\phi$ is an isomorphism. 
Let $y\in X$, and let $\phi^{-1}(y)=y'$.
By Lemma C there exists an apartment $A$ in $X'$ containing $B'$ and $y'$.
The $\phi$-image of $A$ is an apartment in $X$ containing $x$ and $y$.
If two apartments $\sigma(W)$, $\sigma'(W)$ contain $x$ and $y$,
then $(\phi^{-1}\circ\sigma)(W)$ and $(\phi^{-1}\circ\sigma')(W)$
contain $B'$ and $y'$ and thus are isomorphic by an isomorphism $\eta$ fixing 
$B'$ and $y'$. Then $\phi\circ\eta\circ\phi^{-1}$ is an isomorphism 
between $\sigma(W)$, $\sigma'(W)$ fixing $x$ and $y$.
\hfill{$\diamond$(Theorem 4.5)} 
\medskip 
\bf Theorem 4.6\sl\par
For any collection of positive integers $(q_s)_{s\in S}$
there exists a unique $W$-building with $s$-residues of cardinality $q_s+1$.
\rm\smallskip
Proof. 
The construction of the previous subsection with $q_{x,s}=q_s$ for all $x$ 
yields a standard $W$-building $X$ whose residues have required cardinalities.
By Theorem 4.5, $X$ is a building. If $Y$ is another building as in the theorem,
then as in the proof of Theorem 4.2 one can construct 
a morphism $\phi\colon X\to Y$. Moreover, since residue cardinalities agree, 
we can choose $\phi$ so that it is bijective on each $[x]_{\sim_s}$ for
$(x,s)\in R(X)$. By Theorem 4.3 part 4, $\phi$ is a covering map, and hence,
by Lemma A, an isomorphism. 
\hfill{$\diamond$}
\medskip 

One might call the building from Theorem 4.6 a \it regular \rm $W$-building,
with notation $X(W,\q)$ (where $\q=(q_s)_{s\in S}$).
We will now present another construction of $X(W,\q)$.
First, we define an auxiliary local building $Y$.
The set of chambers of $Y$ is a product
$\prod_{s\in S}Y_s$, where each $Y_s$ is a finite set of cardinality at least 2. 
Two chambers $(y_s), (y'_s)$ are $t$-adjacent if $y_s=y'_s$ for all $s\ne t$.
If $T$ is spherical then $Res_Y((y_s),T)$ is isomorphic 
to the product building $\prod_{t\in T}Y_t$ (via dropping coordinates indexed by $S\setminus T$).
Therefore $Y$ is a local $W$-building with the required residue cardinalities,
and $\widetilde{Y}=X(W,\q)$. It follows that $|X(W,\q)|$ carries a free and cocompact action of 
$\Gamma=\pi_1(|Y|)$. If $W$ is a hyperbolic group, then
$|X(W,\q)|$ is $CAT(-1)$ and is quasi-isometric to $\Gamma$; therefore,
$\Gamma$ is Gromov-hyperbolic (in fact, any group acting cocompactly and properly
discontinuously on $|X(W,\q)|$ is Gromov-hyperbolic).

\smallskip

\bf Proposition 4.7\sl\par
The building $|X(W,\q)|$ carries a free and cocompact action of some group $\Gamma$.
If $W$ is hyperbolic, then $\Gamma$ is Gromov-hyperbolic.
\rm\smallskip

Both the proposition and the method of proof (the identification
of $X(W,\q)$ with the universal cover of $Y$) are well-known
(cf. [D2], [GP]).

\medskip

\bf 4.5. Small maps with disjoint images.\rm
\smallskip

A (standard or local) $W$-building is \it thick \rm if each adjacency class has at least 
3 elements.
\smallskip

\bf Theorem 4.8\sl\par
Let $X$ be a thick standard $W$-building, and let $N\subseteq X$ be a finite set. Then 
there exists a finite set $M$, $N\subseteq M\subseteq X$, and two $\pi$-equivariant maps
$\phi,\psi\colon X\to X$ such that $\phi|_M=\psi|_M=Id_M$, $\phi(X\setminus M)\cap\psi(X\setminus M)=\emptyset$.
\rm\smallskip
Proof. Put $M=\pi^{-1}(conv(\pi(N)\cup\{1\}))$.
We claim that $conv(\pi(N)\cup\{1\})$, hence $M$, is finite. Indeed, let $m$ be the number of walls
in $W$ separating some element of $\pi(N)$ from $1$. Let $\ell(w)=p>m$,
and let $w_0=1,w_1,\ldots,w_p=w$ be a minimal gallery. Then one of the $p$ walls between $w_i$ and $w_{i+1}$
does not separate any element of $\pi(N)$ from $1$, while it separates $w$ from $1$; hence,
it separates $w$ from $\pi(N)\cup\{1\}$, and $w\not\in conv(N\cup\{1\})$.
 
Now let $R=R(X)$ be the root set of $X$.
For each $(x,s)\in R$ choose two distinct elements $a_{(x,s)}, b_{(x,s)}\in Res(x,s)\setminus \{x\}$.
Let $A=\{a_r\mid r\in R\}$, $B=\{b_r\mid r\in R\}$.
We now construct $\phi\colon X\to X$ as in the proof of Theorem 4.2. 
Let $x\in X$ be such that $In(\pi(x))=\{s\}$ (i.e., we have a choice for $\phi(x)$).
Let $x_0$ be the shortest element of $Res(x,s)$; then $(x_0,s)\in R$.
If $x\in M$ we put $\phi(x)=x$ (this is allowed, for by induction $\phi(x_0)=x_0$).
If $x\not\in M$, we put $\phi(x)=a_{(\phi(x_0),s)}$.
Similarly we define $\psi$ using $b$'s instead of $a$'s.

Suppose now that $y\in\phi(X)\cap\psi(X)$, $y\not\in M$, and $w=\pi(y)$ is the shortest possible.
Let $y=\phi(x)=\psi(z)$. We have $\phi(x^t)=\psi(z^t)$ for all $t\in In(w)$, hence
$x^t=z^t\in M$ for all such $t$. Now $y\not\in M$ is possible only if $In(w)$ has one element, say $t$.
But then $(x^t,t)\in R$, $(z^t,t)\in R$, $\phi(x)\in A$, $\psi(z)\in B$, contradiction.
\hfill{$\diamond$}

\smallskip
Let $|M|=\{[x,p]\mid x\in M, p\in CL'\}$, $|\phi|([x,p])=[\phi(x),p]$. 
\smallskip
\bf Corollary 4.9\sl\par
Let $\phi,\psi$ be as in Theorem 4.8. Then 
$|\phi|(|X|\setminus |M|)\cap|\psi|(|X|\setminus |M|)=\emptyset$.
\rm\par
Proof. 
Suppose not; let $[x,p]=|\phi|([x_1,p])=|\psi|([x_2,p])$. 
Recall
that $S(p)=\{s\in S\mid \exists \sigma\in L', s\in \sigma, p\in|\sigma|\}$.
Let $y$ be the shortest element of $Res(x,S(p))$, and let $y_1$ be the 
shortest element of $Res(x_1,S(p))$; then $[x,p]=[y,p]$, $[x_1,p]=[y_1,p]$.
Since $Res(x,S(p))$ is the unique residue $R$ of type $S(p)$ in $X$ such that
$|X_{\le R}|$ contains $[x,p]$ (for two different residues $R$ of the same type 
the sets $|X_{\le R}|$ are disjoint), we have $|\phi|(|X_{\le Res(x_1,S(p))}|)
\subseteq |X_{\le Res(x,S(p))}|$, and hence $\phi(Res(x_1, S(p)))\subseteq Res(x,S(p))$.
Now $\pi$-equivariance of $\phi$ implies that $\phi(y_1)=y$. Similarly, 
$\phi(y_2)=y$. It follows that $y\in M$, $[x,p]=[y,p]\in |M|$
\hfill{$\diamond$}
\smallskip

If $X$ is a right-angled hyperbolic building with a folding map $\pi$,
then any $\pi$-equivariant map $\theta\colon X\to X$ fixes the base 
chamber $B=\pi^{-1}(1)$. Therefore $|\theta|\colon|X|\to|X|$ fixes all points in $B$.
Recall that we defined $\partial|X|$ as the space of geodesic rays starting at some base point 
$x_0\in B$. Thus, the map $|\theta|$ induces a continuous map $\partial|\theta|\colon
\partial|X|\to\partial|X|$.

\smallskip
\bf Corollary 4.10\sl\par
Let $X$ be a \ra hyperbolic building, and let $\phi,\psi$ be as in Theorem 4.8.
Then $\partial|\phi|(\partial|X|)\cap\partial|\psi|(\partial|X|)=\emptyset$.
\rm\par
Proof.
Suppose not; let $\partial|X|\ni z=\partial|\phi|(z_1)=\partial|\psi|(z_2)$.
For $y\in \partial|X|$ let $\gamma_y\colon [0,\infty)\to |X|$ be the geodesic from the base point to 
$y$. We have $\gamma_z=|\phi|\circ\gamma_{z_1}=|\psi|\circ\gamma_{z_2}$.
Let $t\in[0,\infty)$ be so large that $\gamma_{z_1}(t)\not\in |M|$ and $\gamma_{z_2}(t)\not\in |M|$;
then $\gamma_z(t)=|\phi|(\gamma_{z_1}(t))=|\psi|(\gamma_{z_2}(t))$,
contradicting Corollary 4.9.
\hfill{$\diamond$}

\jump

\def\x{\overline{x}}
\def\p{\overline{p}}
\def\q{\overline{q}}

\def\ov{\overline}

\centerline{\bf  Appendix}
\medskip
In this appendix we prove an analogue of Lemma 2.15 for arbitrary finite spherical buildings (Theorem A.2).
As a corollary, we deduce that an $n$-dimensional locally finite hyperbolic
or Euclidean building (not necessarily right-angled)
is $(n-2)$-connected at infinity. 
Note that in [GP] 
even more is claimed, but their proof does not convince us. First,
it is not true that $V\cap S(x,s_i+\epsilon)$ 
(here we refer to the proof of Proposition 4.1 in [GP] and we use the notation used 
there)
is of the same homotopy
type as the pointed connected sum of $S(x,s_i-\epsilon)$ with a bouquet of spheres---{\it one for each 
chamber opposite to $c$ in $Lk(y)$}. 
This can be seen by considering 
a $2$-dimensional right-angled building. Second, to claim that $V\cap S(x,s_i+\epsilon)$
has the homotopy type of a bouquet of spheres, one needs to show that $S_+$ is $(n-2)$-connected
(in the notation of [GP]).
This is, in our opinion, a non-trivial fact---Theorem A.2 below. Similar problem appears in
[DM]. Again, as consideration of a $2$-dimensional right-angled building
shows, Lemma 5.5 in [DM] is false. We do not know how to correct this
approach. 

Let $X$ be a finite spherical building of dimension $n\ge1$,
equipped with the standard $CAT(1)$ metric (each apartment is a sphere
of diameter $\pi$). Let $B\in X$ be a chamber, and let $\pi\colon X\to S^n$ 
be the $B$-based folding map. We equip $S^n$ with the standard round metric
such that the restriction of $\pi$ to any apartment $A$ containing $B$ 
isometrically identifies $A$ and $S^n$. The triangulation of $A$ transported by $\pi$ 
is a triangulation of $S^n$; $\pi\colon X\to S^n$ is then a simplicial map.
Images of chambers under $\pi$ will be called chambers.
\medskip
\bf Lemma A.1\sl\par
Let $S_1,S_2,\ldots,S_k\subseteq S^n$ be a finite collection of great 
spheres (of arbitrary dimensions). The set of points $x$ satisfying:

\item{$\bullet$} for every $i$, the function $S_i\ni y\mapsto d(x,y)\in{\bf R}$ has a unique minimum; 
  
\item{$\bullet$} for every $i\ne j$, $d(x,S_i)\ne d(x,S_j)$;
 
is open and dense in $S^n$.\rm
\medskip

For every simplex $\sigma$ in our triangulation of $S^n$ there exists a unique 
smallest great sphere $S\subseteq S^n$ containing $\sigma$.
Apply Lemma A.1 to the collection of all spheres thus obtained; pick a point 
$\x\in S^n$ in the dense open set given by the lemma and inside $int(\pi(B))$.
Let $E^0=\x^\perp\cap S^n$ be the equator of $S^n$ for which $\x$ is a pole;
let $E^-$ be the closed hemisphere with boundary $E^0$ containing $\x$, and let
$E^+$ be the other closed hemisphere with boundary $E^0$. Also, let
$X^0=\pi^{-1}(E^0)$, $X^-=\pi^{-1}(E^-)$, $X^+=\pi^{-1}(E^+)$.
For a subset $U$ or point $p$ of $X$ we put $\overline{U}=\pi(U)$, $\p=\pi(p)$.
We will often pick $\overline{U}$ or $\p$ first, and 
specify $U$ or $p$ later (or at all). For example, $x$ is the unique point in $\pi^{-1}(\x)$.
To avoid a notation clash, closures will be denoted  by $cl$.

\medskip
\bf Theorem A.2\sl\par
Let $X$ be a finite spherical building of dimension $n\ge 1$.
Then $X^+$ is $(n-1)$-connected.\rm
\bigskip 
Proof.
If a set $G\subseteq X$ is isomorphically mapped by $\pi$ onto 
$H\subseteq S^n$, we say that $G$ \it folds onto \rm $H$. 
The following lemma will often be used:
\medskip
\bf Lemma A.3\sl\par
\item{a)} For every chamber $C\in X$ 
there exists an apartment $A$ containing $C$ such that
$A$ folds onto $S^n$.
\item{b)} Let $C_1,C_2$ be two chambers in $X$ such that $\pi(C_2)=-\pi(C_1)$.
Then there exists a unique apartment $A$ containing $C_1$ and $C_2$. This
$A$ folds onto $S^n$.\rm
\smallskip
Proof. a) An apartment containing $C$ and $B$ is good.
b) 
$C_1,C_2$ are opposite in $X$---otherwise $\overline{C_1},\overline{C_2}$ 
would not be opposite in $S^n$. Then $conv(C_1\cup C_2)$ is the desired apartment.
\hfill{$\diamond$}
\smallskip
\bf Remark \rm

One can replace the chamber $C$ in part a) by a point; similarly,
one can replace chambers $C_1,C_2$ in part b) by points $p,q$ such that 
$\q=-\p$ and $p$ is in the interior of some chamber. (Choose chamber/pair of opposite chambers 
containing the point/points, and apply the Lemma.)
\medskip  
\bf Lemma A.4\sl\par
Let $n\ge1$. Then $X^+$ is path-connected.
\rm\smallskip
Proof. Let $p,q\in X^+$. Choose apartments $A_p,A_q$ that fold onto $S^n$ and contain
$p,q$, respectively. Pick a point $p'\in X^0\cap A_p$ which lies in the interior of some chamber $C_p$
(this is possible due to genericity of $\x$). Let $q'\in X^0\cap A_q$
be such that $\overline{q'}=-\overline{p'}$; let $C_q$ be the 
chamber which contains $q$. Then $p$ can be connected to $p'$ by a path in $X^+\cap A_p$, 
and $q$ can be connected to $q'$ by a path in $X^+\cap A_q$. 
Furthermore, $\pi(C_q)=-\pi(C_p)$ so that, by Lemma A.3, there exists an apartment $A\ni C_p,C_q$ which
folds onto $S^n$. Now $p'$ and $q'$ can be connected by a path in $X^+\cap A$. Thus, $p$ and $q$  
can be connected by a path in $X^+$.\hfill{$\diamond$}
\medskip

\bf Proposition A.5\sl\par
Let $n=2$. Then $\pi_1(X)=0$.
\rm\smallskip
Proof. By a general position argument, 
any loop in $X^+$ can be homotoped to a loop           
in $X^+-\pi^{-1}(-\x)$. The latter set deformation retracts onto $X^0$ 
(the deformation retraction moves a point 
along the unique shortest geodesic towards $x$, until it hits $X^0$).
Consequently, any loop in $X^+$ can be homotoped to a loop in $X^0$. 
Now $X^0$ has a natural graph structure, 
inherited form the simplicial structure of $X$. 
Therefore, a loop in $X^0$ is homotopic to a simplicial loop
$\sigma=(e_1,e_2,\ldots,e_k=e_0)$ (each $e_i$ is an oriented 
edge and the endpoint of $e_i$ is the origin of $e_{i+1}$).
A pair $(e_i,e_{i+1})$ will be called a \it backtracking pair \rm (b.p.), if
$\ov{e_i}=\ov{e_{i+1}}^-$ 
(we use $f^-$ to denote $f$ with reversed orientation).  
Now choose an edge $\ov{e}$ in $E^0$. 
A b.p. $(e_i,e_{i+1})$ 
is called \it acceptable, \rm if 
$\ov{e_i}=\pm\ov{e}^\pm$ (one of the four possibilities).
If $(e_i,e_{i+1})$ is a b.p.,
we choose an apartment $A$ that contains $e_i$ and
folds onto $S^n$. Then $A\cap X^0$ is a loop
$(e_i,f_2,f_3,\ldots,f_{2s})$. There exists a smallest
$j\ge 2$ 
such that $\ov{f_j}=\pm\ov{e}^\pm$;
we deform the loop $\sigma=(\ldots,e_i,e_{i+1},\ldots)$
to the loop $\sigma'=(\ldots,e_i,f_2,\ldots,f_j,f_j^-,\ldots,
f_2^-,e_{i+1},\ldots)$. The new loop has the same backtracking pairs
as $\sigma$, with the exception of 
$(e_i,e_{i+1})$, instead of which an acceptable b.p. 
$(f_j,f_j^-)$ appears. 
Notice that the b.p. $(f_j,f_j^-)$ is \it separated, \rm
in the sense that neither $(f_{j-1},f_j)$ nor 
$(f_j^-,f_{j-1}^-)$ is a b.p. 
(if $j=2$, neither $(e_i,f_2)$ nor $(f_2^-,e_{i+1})$ is a b.p.).
Repeating the process, we deform $\sigma$
to a loop with acceptable separated backtracking pairs only.
We keep the notation $\sigma=(e_1,e_2,\ldots,e_k=e_0)$ 
for this new loop.

Now suppose that $\ov{e_i}=\pm\ov{e}^\pm$, but neither
$(e_i,e_{i+1})$ nor $(e_{i-1},e_i)$ is a b.p..
Then $\ov{e_{i+s}}=\ov{e_i}$ so that, by Lemma 2,
there exists an apartment $A\ni e_i,e_{i+s}$. We claim that
$e_{i+1},\ldots,e_{i+s-1}\in A$. 
\smallskip
\bf Lemma A.6\sl\par
Let $\tau=(d_1,\ldots,d_{s+1})$ be a path in $X^0$ such that 
$\ov{d_1}=-\ov{d_{s+1}}$, and let $A$ be the apartment 
in $X$ containing $d_1$ and $d_{s+1}$. Then $\tau$ is contained in 
$A$.\par\rm
Proof. The path
$\tau'=(d_2,\ldots,d_s)$ from the endpoint $y$ of $d_1$ 
to the origin $z$ of $d_{s+1}$ has geometric length
$d(\ov{y},\ov{z})$. Since $d(\ov{y},\ov{z})\le d(y,z)$,
$\tau'$ is a shortest geodesic. Now $y,z\in A$, and $A$ is convex,
therefore $\tau$ is contained in $A$. \hfill{$\diamond$(Lemma A.6)}
\smallskip

In $A\cap X^+$, the path $(e_{i+1},e_{i+2},\ldots,e_{i+s-1})$ is
homotopic (with endpoints held fixed) to a path 
$(e_i^-,f_{i+1},\ldots,f_{i+s-1},e_{i+s}^-)$, where
$\ov{f_{i+j}}=-\ov{e_{i+s-j}}^-$. The effect of this change 
on $\sigma$ is $$(\ldots,e_i,e_{i+1},\ldots,e_{i+s-1},e_{i+s},\ldots)
\rightarrow
(\ldots,e_i,e_i^-,f_{i+2},\ldots,f_{i+s-1},e_{i+s}^-,e_{i+s},\ldots).$$
It may happen that $(e_{i+s},e_{i+s+1})$ is a b.p.; if this is the case,
we further modify the loop: 
$$(\ldots,f_{i+s-1},e_{i+s}^-,e_{i+s},e_{i+s+1},\ldots)
\rightarrow (\ldots,f_{i+s-1},e_{i+s+1},\ldots).$$  
Travelling along the loop and repeating the process if necessary,
we finally arrive at a loop
$\sigma=\sigma_1\sigma_2\ldots\sigma_{2l}$, where
each $\sigma_i$ is a path of length $s+1$ with no b.p.,
and (last edge of $\sigma_i$, first edge of $\sigma_{i+1}$)
is an acceptable separated b.p. (for $i=0,1,\ldots,2l-1$, where
$\sigma_0=\sigma_{2l}$). 

Suppose now that $\tau$ is a path of length $s+1$ containing no 
b.p.. Let $A$ be the apartment containing the extreme edges of $\tau$.
Then $\tau$ is homotopic (with the endpoints held fixed)
in $A\cap X^+$ to a path
$\hat{\tau}\subseteq A\cap X^0$ of length $s-1$.
Now we modify $\sigma$ by homotopy inside $X^+$,
changing $\sigma_j$ and $\sigma_{2l-j+1}$ to
$\hat{\sigma}_j$, $\hat{\sigma}_{2l-j+1}$ (resp.), for 
all positive even $j\le l$. We obtain a loop $\sigma=\eta\xi$ with exactly two 
backtracking pairs, where $\eta$, $\xi$ are paths of equal length,
say $u$,
and none of them  contains a b.p.. A loop of this form 
will be called a \it $u$-moon.\rm
\smallskip
\bf Lemma A.7\sl\par
An $(s+1)$-moon is contractible in $X^+$.\rm\par
Proof. If $u=s+1$, then $\eta\xi$ is homotopic (in $X^+$) 
to $\eta\hat{\xi}$; the latter is contained, by Lemma 4, in 
the apartment $A$ spanned by the extreme edges of $\eta$.
The apartment $A$ folds onto $S^n$, therefore 
$\eta\hat{\xi}$ is null-homotopic in $A\cap X^+$.
\hfill\hbox{$\diamond$(Lemma A.7)}
\smallskip
\bf Lemma A.8\sl\par
If $u>s+1$, then a $u$-moon is homotopic to
a concatenation of an $(s+1)$-moon and a $(u-1)$-moon.\rm\par
Proof.
Let $\eta=(\eta_1,\ldots,\eta_u)$,
$\xi=(\xi_1,\ldots,\xi_u)$. Let $A$ be the apartment 
spanned by $\eta_1$ and $\xi_{u-s}$, and let
$\tau$ be the path of length $s-1$ in $A\cap X^0$ from
the endpoint of $\eta_1$ to the endpoint of $\xi_{u-s}$.
Then $\eta\xi$ is homotopic to the concatenation
of the $(s+1)$-moon  $\tau\xi_{u-s}^-\xi_{u-s}\xi_{u-s+1}\ldots\xi_u\eta_1$
and the $(u-1)$-moon $\eta_2\eta_3\ldots\eta_u\xi_1\xi_2\ldots\xi_{u-s}
\tau^-$. \hfill\hbox{$\diamond$(Lemma A.8)}
\smallskip

Repeated application of Lemmas A.8 and A.7 finishes 
the proof of Proposition A.5.
\null\nobreak\hfill{$\diamond$(Proposition A.5)}
\medskip

Thus, Theorem A.2 is true for $n=1,2$. We will proceed with the proof of the general 
case by induction on $n$. Suppose that $n>2$, and that Theorem A.2 is true for
all finite buildings of dimension less than $n$. Let $X$ be a finite building of dimension $n$.

Let $\sigma_1,\ldots,\sigma_\ell$ be all the  simplices of our triangulation of $S^n$ 
that have the following property:
there exists a minimal unit-speed geodesic $\ov{\gamma_i}$ issued from
$\ov{x}$ which intersects the interior of $\sigma_i$ at $\ov{p_i}=\ov{\gamma_i}(t_i)$.
By the choice of $\ov{x}$, $\ov{\gamma_i}$ is unique, all the $t_i$ are distinct and 
none of them equals $\pi/2$. We can assume that $t_1<\ldots<t_\ell$.
Let $X^+_r=\{y\in X\mid d(\pi(y),\x)\ge r\}$. Our strategy is to show, by induction 
on $i$, that $X^+_r$ is $(n-1)$-connected for $t_i<r<t_{i+1}$, $r\le\pi/2$.
To this end, we need to prove that $X^+_\epsilon$ is $(n-1)$-connected for 
small positive $\epsilon$, and then we need to understand how $X^+_r$ changes when
$r$ switches from the interval $(t_{i-1},t_i)$ to $(t_i,t_{i+1})$.

If $\epsilon$ is sufficiently close to $0$, then $X^+_\epsilon$ is homotopy
equivalent to $X-\{x\}$, which homotopically is a bouquet of $n$-spheres with one sphere
punctured, and so is $(n-1)$-connected.

Now we will closely follow the proof of Lemma 3.3.  
Put $t=t_i$, $\sigma=\sigma_i$, $\ov{\gamma}=\ov{\gamma_i}$, $\ov{p}=\ov{p_i}$. 
We choose a small $\delta>0$ such that the sphere $S_\delta(\ov{p})$ is contained in $int(Res\,\sigma)$, 
where $Res\,\sigma=\bigcup\{\tau\mid\sigma\subseteq\tau\}$.
Then there exists an $\epsilon>0$ such that
$S_{t+\epsilon}(\x)\cap\sigma\subseteq S_\delta(\ov{p})$.
The intersection $S_{t+\epsilon}(\x)\cap\sigma$ is a sphere of dimension $d=\dim{\sigma}-1$
contained in the interior of $\sigma$.
Decreasing $\delta$ we ensure that the following inequalities hold: 
$\epsilon<\min\{t_{i+1}-t, t-t_{i-1}\}$, $t+\epsilon<\pi/2$.
Let $\ov{D}$ denote $S_\delta(\ov{p})\cap\sigma^\perp_{\ov{p}}$
(where by $\eta^\perp_{\ov{p}}$ we denote the largest great sphere in $S^n$ orthogonal to $\eta$ at $\ov{p}$).
The intersection $\ov{D}\cap\ov{\gamma}^\perp_{\ov{p}}$ 
divides $\ov{D}$ into two closed hemispheres:
$\ov{D}^-$ (the one closer to $\x$) and $\ov{D}^+$.
Observe that $\ov{D}^+=\ov{D}\setminus B_{t+\epsilon}(\x)$.
The sphere $\ov{D}$ inherits a triangulation from $S^n$. We want 
the spheres $S_{t\pm\epsilon}(\ov{x})$ to intersect this triangulation
`in the same way'. This condition can be achieved by further decreasing $\delta$ 
(and, consequently, $\epsilon$). 
\smallskip
 
Now we pass to $X$. Let 
$$\eqalignno{
\pi^{-1}(\sigma)&=\{\sigma_1,\ldots,\sigma_k\},\cr
\{p_s\}&=\sigma_s\cap\pi^{-1}(\p),\cr
K_s&=\{y\in X\mid d(y,p_s)\le\delta\},\cr
K&=\coprod_sK_s,\cr 
Y^+&=cl(X^+_{t+\epsilon}\setminus K),\cr
D_s&=\pi^{-1}(\ov{D})\cap K_s,\cr
D_s^+&=D_s\cap X^+_{t+\epsilon},\cr
S_s^d&=\sigma_s\cap S_\delta(p_s).}$$
The first two definitions override our previous convention. 
Observe that: $X^+_{t+\epsilon}$ is homotopy equivalent 
to $Y^+$; $X^+_{t-\epsilon}$ is homotopy equivalent to $Y^+\cup K$;
$K_s$ is homeomorphic to a cone over the join
$S^d_s*D_s$, and is attached to $Y^+$ along a subset
homeomorphic to $S^d_s*D^+_s$.

Each $D_s$ is a spherical building of dimension $n-d-2$.
By the inductive assumption $D^+_s$ is $(n-d-3)$-connected, 
which implies that $S^d_s*D^+_s$
is $(n-2)$-connected. Moreover, $K_s$ is contractible. 
Van Kampen's theorem, Mayer--Vietoris sequence and the inductive 
assumption that $Y^+\cup K$ is $(n-1)$-connected imply 
$(n-2)$-connectedness of $Y^+$. It remains to prove that 
$H_{n-1}(Y^+)=0$ 
(here we depart from the proof of Lemma 3.3).
It follows from the Mayer-Vietoris sequence that 
$H_{n-1}(Y^+)$ is generated by the images of 
$H_{n-1}(S^d_s*D^+_s)$ ($s=1,\ldots,k$). We will show that any $(n-1)$-cycle 
in $S^d_s*D^+_s$ 
is null-homologous in $Y^+$. 

Let us subdivide the usual triangulation of $S^n$ to a minimal
cellulation in which $S_{t+\epsilon}(\x)\setminus B_\delta(\ov{p})$ and $S_\delta(\ov{p})$ are subcomplexes.
Pull this cellulation back to $X$ via $\pi$.
An $(n-1)$-cycle $z$ in $S^d_s*D^+_s$ is a join of 
the fundamental class of $S^d_s$ and an $(n-d-2)$-cycle $\tilde{z}$ in $D^+_s$.
The cycle $\tilde{z}$ can be regarded as a cycle in $D_s$ vanishing 
outside $D^+_s$. Now every cycle in $D_s$ can be expressed as a combination 
of fundamental classes of apartments. More specifically,
let $c_0$ be the chamber in $D_s$ which is closest to $\x$, and let 
$c_1,\ldots,c_m$ be all the chambers in $D_s$ opposite to $c_0$.
Let $a_i$ be the apartment in $D_s$ containing $c_0$ and $c_i$;
then $\tilde{z}=\sum_{i=1}^m\alpha_i[a_i]$ for some integers $\alpha_i$.
Let $C_i$ be the chamber of $X$ containing $c_i$. Let
$C_{-1}$ be a chamber in $X$ such that $\ov{C_{-1}}=-\ov{C_1}$ 
(note that $\ov{C_1}=\ldots=\ov{C_m}$),
and let $A_i$ be the apartment in $X$ containing $C_{-1}$ and $C_i$.
Let $Z=\sum_{i=0}^m\alpha_i[A_i]\in Z_n(X)$.
Split $Z$ into $Z_1+Z_2$, where $Z_1\in C_n(K_s)$, 
$Z_2\in C_n(cl(X\setminus K_s))$.    
Clearly, $Z_1$ is the cone over $z$ so that $\partial Z_1=z$.
Therefore, $\partial(-Z_2)=z$. We claim that $Z_2\in C_n(Y^+)$.
First, notice that if $u\ne s$ then for all $i$ we have $K_u\cap A_i=\emptyset$:
since $A_i$ folds onto $S^n$, it can intersect only one component of
$K$, and it does intersect $K_s$.  
Next, let $C$ be a chamber of $X$ contained in $A_i$ but not contained in 
$Y^+\cup K$. Choose a point $y\in int(C)$ such that $y\not\in X^+_t\cup K$
(in particular, $\overline{y}\ne-\p$). The unique minimal geodesic
$\gamma$ from $p_s$ to $y$ is contained in
$(A_i\cap(X\setminus X^+_t))\cup\{p_s\}$, because both extremities 
belong to this convex set.
In particular, $\gamma$ intersects $S_\delta(p_s)$ outside
$X_t^+$, hence outside $Y^+$. It follows that $\gamma$ 
leaves $p_s$ through the interior of a chamber $C'$
on which $Z$ is zero. The chamber $C'$ is the closest
to $C$ (in the gallery distance) among all chambers
in the residue of $\sigma$. Now all apartments $A_j$ contain
$\sigma$; therefore, if an apartment 
$A_j$ contains $C$, it also contains $C'$.
Conversely, we claim that if $C'\in A_j$, then $C\in A_j$.
To see this extend the geodesic $\gamma$ to
$\gamma\colon[0,\pi]\to A_i$ (so that $\gamma(\pi)$ is opposite
to $\gamma(0)=p_s$ in $A_i$).
Now slightly rotate $\gamma$ inside $A_i$, so as to obtain
a geodesic $\eta$ which still passes through the interiors of
$C'$ and $C$ ($\eta(a)\in int(C'),\eta(b)\in int(C), a<b$),
but $\eta(0)\in int(C_i),\eta(\pi)\in int(C_{-1})$.
Then $\eta|_{[a,\pi]}$ is a minimal geodesic.
Suppose that  $C'\in A_j$. Then $\eta(a),\eta(\pi)\in A_j$
so that $\eta(b)\in A_j$, hence $C\in A_j$.
Thus, we have verified that
$\{j\mid C\in A_j\}=\{j\mid C'\in A_j\}$. The value of $Z$ on $C$ is equal to
$\sum_{j\mid C\in A_j}\alpha_j=\sum_{j\mid C'\in A_j}\alpha_j$; the latter
is the value of $Z$ on $C'$, i.e. zero.
\hfill{$\diamond$(Theorem A.2)}   
\medskip

A space $X$ is \it $k$-connected at infinity, \rm if for every compact $K\subseteq X$ 
there exists a compact $L$,  $K\subseteq L\subseteq X$, such that any map 
$S^i\to X\setminus L$ extends to a map $B^{i+1}\to X\setminus K$ 
(for $i=0,1,\ldots,k$). 
 
\bf Corollary A.9\sl\par
An $n$-dimensional locally finite hyperbolic or Euclidean 
building is $(n-2)$-connected
at infinity.
\rm\par
Proof. It is enough to check that complements of balls $B_r$ are $(n-2)$-connected, for 
$r>0$ arbitrarily large.
By the geodesic retraction such a complement is homotopically equivalent
to $S_r$. Then the proof of Lemma 3.3, for $U=X\cup \partial X$, goes through, with Theorem A.2 used instead
of Lemma 2.15. Lemma 2.16 is never needed for this choice of $U$.
\hfill{$\diamond$} 

\jump

\centerline{\bf References}
\medskip

\item{[B]} Benakli N., Poly\` edres hyperboliques, passage du local au global, Th\` ese, 
Universit\' e Paris Sud, 1992.  

\item{[BK]} Benakli N., Kapovich I.,
Boundaries of hyperbolic groups,
Combinatorial and geometric group theory (New York, 2000/Hoboken, NJ, 2001),  
39--93, Contemp. Math., 296, Amer. Math. Soc., Providence, RI, 2002. 

\item{[Be]} Bestvina M., Characterizing
$k$-dimensional universal Menger compacta, Mem. Amer. Math. Soc. 71
(1988), no. 380.

\item{[Bb]} Bourbaki N., 
``Groupes et algebres de Lie, chapitres IV-VI", Hermann, 1968.

\item{[BH]} Bridson M., Haefliger A., 
``Metric spaces of non-positive curvature'', Springer, 1999.

\item{[Bd1]}
Bourdon M., Immeubles hyperboliques, dimension conforme et rigidit\' e 
de Mostow, Geom. Funct. Anal. {\bf 7} (1997), no. 2, 245--268. 

\item{[Bd2]}
Bourdon M., 
Sur les immeubles fuchsiens et leur type de quasi--isom\'etrie, Ergodic Theory 
and Dynamical Systems
{\bf 20} (2000), no. 2, 343--364.


\item{[BP1]} Bourdon M., Pajot H., Rigidity of
quasi--isometries for some hyperbolic buildings, 
Comment. Math. Helv. {\bf 75} (2000), no. 4, 701--736.

\item{[BP2]} Bourdon M., Pajot H., Cohomologie $L^p$ et \' espaces de Besov, 
J. Reine Angew. Math. {\bf 558} (2003), 85--108.

\item{[BMcCM]}
Brady N., McCammond J., Meier.J., 
Local--to--asymptotic topology for cocompact $\rm CAT(0)$ complexes,  
Topology Appl.  {\bf 131}  (2003),  no. 2, 177--188.

\item{[Br]} Brown K., Buildings, Springer, New York, 1989. 

\item{[D1]} Davis M.W., Groups generated by reflections 
and aspherical manifolds not covered by Euclidean space,  
Ann. of Math. (2)  {\bf 117}  (1983),  no. 2, 293--324. 

\item{[D2]} Davis M.W., 
Buildings are CAT(0), 
in "Geometry and Cohomology in Group Theory", 
edited by P. Kropholler, G. Niblo, R. Stohr,
London Math. Society Lecture Notes \bf 252\rm, 
Cambridge Univ. Press, Cambridge, 1998, pp.108--123.


\item{[DM]} Davis M.W., Meier, J., 
The topology at infinity of Coxeter groups and buildings,  
Comment. Math. Helv.  {\bf 77} (2002),  no. 4, 746--766. 

\item{[Gl]} Globus M., unpublished, 1996.

\item{[GP]}
Gaboriau D., Paulin F., Sur
les immeubles hyperboliques,
Geom. Dedicata {\bf 88} (2001), no. 1-3, 153--197.

\item{[HP]}
Haglund F., Paulin F., 
Constructions arborescentes d'immeubles,
Math. Ann. {\bf 325} (2003), no. 1, 137--164.

\item{[KK]}
Kapovich, M., Kleiner, B.,
Hyperbolic groups with low-dimensional boundary, 
Ann. Sci. \' Ecole Norm. Sup. (4) {\bf 33} (2000), no. 5, 647--669.

\item{[M]} Moussong G., Hyperbolic Coxeter groups, PhD thesis, 
the Ohio State University, 1987. 

\item{[Ron]} Ronan M.,
Lectures on buildings, Perspectives in Mathematics vol. 7,
Academic Press 1989. 

\item{[Se]} Serre J.-P., Cohomologie des groupes discrets,
in Prospects in mathematics, 
Ann. of Math. Studies, No. 70, Princeton University Press, Princeton, N.J.,
 1971, pp.77--169.

\item{[Vin]} Vinberg E.B., Hyperbolic reflection groups,
Usp. Mat. Nauk {\bf 40} (1) (1985), 29--66.
\medskip

\bye